\documentclass[a4paper,10pt]{amsart}
\usepackage[arrow,matrix]{xy}
\usepackage{amsmath,amssymb,amscd,bbm,amsthm,mathrsfs}

\theoremstyle{plain}
\newtheorem{thm}{Theorem}[subsection]

\newtheorem{lem}[thm]{Lemma}
\newtheorem{prop}[thm]{Proposition}

\newtheorem{exa}[thm]{Example}

\newtheorem{note}[thm]{Notations}

\newtheorem{rem}[thm]{Remark}

   \def\op{\oplus} \def\ot{\otimes}
\def\Hom{\operatorname {Hom}}

\def\Ext{\operatorname {Ext}} 
 \def\k{\mathbbm{k}}

\def\id{\operatorname{id}}

\begin{document}
\title{\bf PBW deformations of Koszul algebras over a nonsemisimple ring}

\author{Ji-Wei He, Fred Van Oystaeyen and Yinhuo Zhang}
\address{J.-W. He\newline \indent Department of Mathematics, Shaoxing College of Arts and Sciences, Shaoxing, Zhejiang 312000,
China} \email{jwhe@usx.edu.cn}
\address{F. Van Oystaeyen\newline\indent Department of Mathematics and Computer
Science, University of Antwerp, Middelheimlaan 1, B-2020 Antwerp,
Belgium} \email{fred.vanoystaeyen@ua.ac.be}
\address{Y. Zhang\newline
\indent Department WNI, University of Hasselt, Universitaire Campus,
3590 Diepenbeek, Belgium} \email{yinhuo.zhang@uhasselt.be}

\date{}
\begin{abstract}
Let $B$ be a generalized Koszul algebra over a finite dimensional algebra $S$. We construct a bimodule Koszul resolution of $B$ when the projective dimension of $S_B$ equals 2. Using this we prove a Poincar\'{e}-Birkhoff-Witt (PBW) type theorem for a deformation of a generalized Koszul algebra. When the projective dimension of $S_B$ is greater than 2, we construct bimodule Koszul resolutions for generalized smash product algebras obtained from braidings between finite dimensional algebras and Koszul algebras, and then prove the PBW type theorem. The results obtained can be applied to standard Koszul Artin-Schelter Gorenstein algebras in the sense of Minamoto and Mori \cite{MM}.
\end{abstract}


\maketitle

\section*{Introduction}

Let $\k$ be a field of characteristic $p$ and let $G\subseteq SL(2,\k)$ be a cyclic group of order $p$. $G$ acts on the polynomial algebra $\k[x,y]$. Let $\k G$ be the group ring. Denote by $\k[x,y]\#\k G$ the skew group algebra. For an element $\lambda\in \k G$, set $H_\lambda=\k\langle x,y\rangle/(xy-yx-\lambda)$, where $\k\langle x,y\rangle$ is the free algebra generated by $x,y$. The algebra $H_\lambda$ is called a symplectic reflection algebra in positive characteristic. Norton proved in \cite{N} that $H_\lambda$ is a Poincar\'{e}-Birkhoff-Witt (PBW) deformation of the skew group algebra $\k[x,y]\#\k G$. This does not occur coincidently. Shepler and Witherspoon proved in \cite{SW} (more recently \cite{SW2}) the PBW theorem for any skew group algebra obtained from an action of a finite group on a Koszul algebra. More generally, Walton and Witherspoon \cite{WW} proved the PBW theorem for any smash product algebra obtained from a Hopf algebra acting on a Koszul algebra.

Note that smash product algebras stemming from Hopf actions on Koszul algebras are {\it generalized Koszul} (for the definition, see Section \ref{sec0}) in the sense of Woodcock \cite{W}, which are Koszul analogous algebras defined over an arbitrary ring. We are wondering whether the PBW theorem still holds for generalized Koszul algebras. The PBW theory for the classical Koszul algebras was developed by many authors, see, e.g. \cite{BG,PP}. As we know, PBW conditions are determined by the Hochschild cohomology, and hence the crucial point in the proof of the PBW theorem for a Koszul algebra is the construction of the bimodule Koszul resolution. Unlike the classical case, the base ring of a generalized Koszul algebra is no longer semisimple. The cohomolgy of the base ring will complicate the cohomology of the algebra. Hence it is usually complicated and difficult to construct the bimodule Koszul resolution of a generalized Koszul algebra, which in turn makes the proof of the PBW theorem for a generalized Koszul algebra much more difficult.

However, if the projective dimension of the base ring, viewed as a graded module over the generalized Koszul algebra, is not larger than 2, then, as we will see in Section \ref{sec0}, the cohomology of the base ring will not complicate the cohomology of the generalized Koszul algebra too horribly. In this case, we can construct a bimodule projective resolution for the generalized Koszul algebra, and are able to prove the PBW theorem.

Let $S$ be a finite dimensional  algebra, ${}_SM_S$ a finite dimensional  $S$-bimodule. Let $T_S(M)$ be the tensor algebra over $S$. Assume that $B=T_S(M)/(R)$ is a generalized Koszul algebra where $R\subseteq M\ot_S M$ is a sub-$S$-bimodule. Take a linear map $\theta:R\to S$. Define $U=T_S(M)/(r-\theta(r):r\in R)$. Then $U$ is a filtered algebra with the natural ascending filtration. If the projective dimension of $S_B$ is 2, we provide a method for constructing the bimodule Koszul resolution of $B$ in Section \ref{sec}. We then prove the following result (cf. Theorem \ref{thm3}).

\noindent{\bf Theorem A.} {\it  $U$ is a PBW deformation of $B$ if and only if $\theta$ is an $S$-bimodule morphism.}

Note that in the theorem above the generating relations of the algebra $U$ do not involve the degree 1 part of the tensor algebra. We remind the reader that the PBW deformations appeared in \cite{BG,PP,SW,SW2,WW} are more general in terms of generating relations. However, the situation in the theorem above occurs very common, and it applies to many interesting algebras, such as symplectic reflection algebras (cf. \cite{CBH,EG,N}), graded Hecke algebras (cf. \cite{EG,RS,NW}) for arbitrary finite group rings, and Lusztig-type algebras in positive characteristic (cf. \cite{Lu}, also \cite{SW2}),...,etc. Our result provides an alternative approach to those algebras.

If the projective dimension of $S_B$ is greater than 2, then the method used in Section \ref{sec} no longer works. In what follows, we write  $B=\op_{i\ge0} B_i$ for a generalized Koszul algebra. In order to construct the bimodule Koszul resolution for a Koszul algebra $B$ with pdim$S_B>2$, we assume further that $B_i$ is free both as a left $S$-module and as a right $S$-module for all $i\ge1$. Note that $B$ is generated by $M$ over $S$. Since $M$ is free as a right $S$-module, we may assume $M=V\ot S$ as a free right $S$-module for some finite dimensional vector space $V$. The left action of $S$ on $M$ induces a braiding $\Psi:S\ot V\to V\ot S$ (for details, see Subsection \ref{sec1}). This yields a generalized smash product $T(V)\#S$ such that $B$ is a quotient algebra of $T(V)\#S$. Hence in Section \ref{secc}, we consider the PBW deformations of a generalized smash product algebra obtained from a braiding between $S$ and a classical Koszul algebra $A$. We provide a method for constructing the bimodule Koszul resolution of the generalized smash product algebra $A\#S$ in Subsection \ref{sec3} (cf. Theorem \ref{thm}), and then prove the following result (cf. Theorem \ref{thm2}, for the notions see Subsection \ref{sec4}), generalizing the corresponding results in \cite{GGV,SW,SW2,WW}.

\noindent{\bf Theorem B.} {\rm Let $U=T_S(M)/(r-\phi_B(r)-\theta_B(r):r\in R)$. Then $U$ is a PBW deformation of $B=A\#S$ if and only if the following conditions are satisfied:

{\rm(i)} $(\phi_B\ot_S \id-\id\ot_S\phi_B)(\overline{R}\ot_S M\cap M\ot_S\overline{R})\subseteq \overline{R}$;

{\rm(ii)} $\phi_B(\phi_B\ot_S \id-\id\ot_S\phi_B)=-(\theta_B\ot_S \id-\id\ot_S \theta_B)$;

{\rm(iii)} $\theta_B(\id\ot_S\phi_B-\phi_B\ot \id)=0$.}

In \cite{MM} Minamoto and Mori introduced the notion of Artin-Schelter Gorenstein algebra over nonsemisimple algebras. We study the PBW deformation of Koszul Artin-Schelter Gorenstein algebras in Section \ref{sec5}. We show that if $B$ is a standard Artin-Schelter Gorenstein algebra, then $B_0$ is selfinjective and the projective dimension of $B_0$ viewed as a graded $B$-module is finite (cf. Proposition \ref{prop10}). This generalizes a similar result in \cite{DW} for classical Koszul algebras. As the main result of Section \ref{sec5}, we prove the PBW theorem for a standard Koszul Artin-Schelter Gorenstein algebra of graded injective dimension 2 (cf. Theorem \ref{thm4}).

Throughout the paper, $\k$ is a fixed field. Unadorned $\ot$ means $\ot_\k$, and $\Hom$ means $\Hom_\k$. All the vector spaces, algebras and modules are defined over $\k$.

\section{Generalized Koszul algebras}\label{sec0}

Let $B=B_0\op B_1\op\cdots$ be a positively graded algebra assumed to be locally finite, that is, $\dim B_i<\infty$ for all $i\ge0$. The definition of the classic Koszulity of $B$ requires  $B_0$ to be semisimple \cite{P,BGS}. This requirement is sometimes too strong. Many interesting algebras with properties similar to a Koszul algebra do not satisfy this requirement. For example, the recent study of the structure of Artin-Schelter regular algebras requires $B_0$ not to be a semisimple algebra \cite{MM}. A generalization of the Koszulity is necessary. In fact, there exist several versions of a generalized Koszulity, see \cite{W,GRS,L,M}. In this paper, we use Woodcock's definition \cite{W}, that is, a positively graded algebra $B$ is a {\it right Koszul algebra} if $B_i$ is projective both as a left $S$-module and as a right $S$-module for all $i\ge1$, and $B_0$ viewed as a right graded $B$-module has a graded projective resolution
$$\cdots \longrightarrow P^{-n}\longrightarrow\cdots \longrightarrow P^{-1}\longrightarrow P^0\longrightarrow B_0\longrightarrow 0,$$ such that $P_n$ is generated in degree $n$ for all $n\ge0$. If $B_0=\k$, we say that $B$ is {\it classically Koszul}. A {\it left Koszul algebra} is defined similarly. Woodcock proved in \cite{W} that $B$ is left Koszul if and only if it is right Koszul. Henceforth, we omit the prefix ``left'' or ``right''.

\subsection{Resolutions of Koszul algebras}
Let $B=B_0\op B_1\op\cdots$ be a Koszul algebra. Set $S=B_0$, and $M=B_1$. Let $T_S(M)=S\op M\op M\ot_S M\op\cdots$ be the tensor algebra over $S$. $B$ is a quadratic algebra in the sense that $B= T_S(M)/I$ where $I$ is a two-sided ideal generated by a set of elements in $M\ot_S M$ \cite{W,L}. Let $\pi:T_S(M)\to B$ be the projection map. Let $R$ be the kernel of the restriction of the projection $\pi$ to $M\ot_SM\to B_2$. Then $R$ is an $S$-bimodule, and the two-sided ideal $I$ is generated by $R$.

Since $B$ is projective as an $S$-module on both sides, we may view $M\ot_S R$ and $R\ot_SM$ as submodules of $M\ot_S M\ot_SM$ through the natural inclusion map. Therefore, $M\ot_S R\cap R\ot_SM$ makes sense. In general, for $i,j\ge0$, $M^{\ot^i_S}\ot_S R\ot_S M^{\ot^j_S}$ is a submodule of $M^{\ot_S^{i+j+2}}$, and $\bigcap_{i=0}^{n-2}M^{\ot_S^ i}\ot_S R\ot_S M^{\ot_S^{n-2-i}}$ is a sub-bimodule of $M^{\ot_S^n}$ for all $n\ge2$. A Koszul resolution for the graded right $B$-module $S$ was constructed in \cite{W}. The following result was proved in \cite{W}. We include our own proof for later use.

\begin{prop}\label{prop7} \cite[Theorem 2.3]{W} Let $B$ be a Koszul algebra. Set $M=B_1$ and $S=B_0$.
We have the following graded projective resolution of the graded module $S_B$:
$$\cdots \longrightarrow K_n\ot_S B \overset{\partial^{-n}}\longrightarrow\cdots\longrightarrow K_2\ot_S B\overset{\partial^{-2}}\longrightarrow K_1\ot_SB\overset{\partial^{-1}}\longrightarrow B\overset{\varepsilon}\longrightarrow S_B\longrightarrow0,$$ where $$K_n=\bigcap_{i=0}^{n-2}M^{\ot_S^i}\ot_S R\ot_S M^{\ot_S^{n-2-i}}$$ for all $n\ge2$ and $K_1=M$, and the differential $\partial^{-n}$ is the restriction of the right $B$-module morphism $M^{\ot_S^n}\ot_SB\longrightarrow M^{\ot_S^{ n-1}}\ot_SB$ defined by $x_1\ot_S\cdots\ot_S x_n\ot_S b\mapsto x_1\ot_S\cdots\ot_S x_{n-1}\ot_S x_nb$ for $x_1,\dots,x_n\in M,b\in B$.
\end{prop}
\proof Note that the consequence $M\ot_SB\overset{\partial^{-1}}\longrightarrow B\overset{\varepsilon}\longrightarrow S\longrightarrow0$ is clearly exact. Suppose that we already have an exact sequence $K_n\ot_S B \overset{\partial^{-n}}\longrightarrow\cdots\longrightarrow K_1\ot_SB\overset{\partial^{-1}}\longrightarrow B\overset{\varepsilon}\longrightarrow S_B\longrightarrow0$, where $K_i\ot_SB$ is graded projective for all $0\leq i\leq n$. By the definition of the differential $\partial^{-n}$, the degree $n+1$ part of ker$\partial^{-n}$ is $K_n\ot_SM\cap M^{\ot_S^{n-1}}\ot_S R$ which equals  $K_{n+1}$. Since $B$ is Koszul and $K_{i}\ot_SB$ is a graded projective cover of the kernel of $\partial^{-i+1}$ for all $n-1\ge i\ge0$, we see that ker$\partial^{-n}$ is generated in degree $n+1$. Since $B$ is projective as a right $S$-module, the degree $n+1$ part of the sequence
\begin{equation}\label{eq23}
  0\longrightarrow\text{ker}\partial^{-n}\longrightarrow K_n\ot_S B \overset{\partial^{-n}}\longrightarrow\cdots\longrightarrow K_1\ot_SB\overset{\partial^{-1}}\longrightarrow B\overset{\varepsilon}\longrightarrow S_B\longrightarrow0
\end{equation}
is split exact. Hence $K_{n+1}$ is projective as a right $S$-module. Therefore, $K_{n+1}\ot_SB\longrightarrow \text{ker}\partial^{-n}$ is a graded projective cover. \qed

\begin{rem}\label{rem1}{\rm Note that $K_n$ is an $S$-bimodule for all $n\ge0$. It is not hard to see that the differentials in the Koszul resolution constructed above is also compatible with the left $S$-module structures. By the exact sequence (\ref{eq23}) as in the proof of Proposition \ref{prop7}, one can show inductively that $K_n$ is also projective as a left $S$-module since $B$ is also projective as a left $S$-module.}
\end{rem}

\subsection{Bimodule resolutions}
From the Koszul resolution constructed above, we may define a $B$-bimodule (not necessary projective) resolution of $B$. We need the following lemma.

\begin{lem}\label{lem5} Let $B$ be a positively graded algebra with $B_0=S$ such that $B$ is projective as a right $S$-module. Let $\{V^{i}\}$ ($i\geq 0$) be a sequence of left graded $S$-modules which are concentrated in nonnegative degrees. Let $\cdots\longrightarrow B\ot_SV^i\longrightarrow\cdots\longrightarrow B\ot_S V^1\longrightarrow B\ot_S V^0\longrightarrow0$ be a complex of graded left $B$-modules. If the resulting complex after applying the functor $S\ot_B-$ on the sequence is exact, then the sequence itself is exact.
\end{lem}
\proof Denote by $W^{-i}$ the left $B$-module $B\ot_S V^i$ and by $W^\bullet$  the complex in the lemma. We define a filtration on $W^\bullet$ as follows. By assumption, $V^i=\op_{j\ge0}V^i_j$ for all $i\ge0$. For $i\ge0$, set $F_jW^i=0$ for all $j<0$, $F_0W^{-i}=B\ot_S V^i_0$, $F_1W^{-i}=B\ot_S(V^i_0\op V^i_1),\dots, F_nW^\bullet=B\ot_S (V^i_0\op\cdots\op V^i_n),\dots$. Then it is easy to see that $F_nW^\bullet$ is a subcomplex of $W^\bullet$ for all $n\ge0$, moreover, $\bigcup_{n\ge0}F_nW^\bullet=W^\bullet$. Applying the functor $S\ot_B-$ on the complex $W^\bullet$, we obtain an exact sequence of graded $S$-modules $$\cdots\longrightarrow V^i\overset{f^i}\longrightarrow\cdots \longrightarrow V^1\overset{f^1}\longrightarrow V^0\longrightarrow0.$$ Hence for each $j\ge0$, the sequence $$\cdots\longrightarrow V^i_j\overset{f^i}\longrightarrow\cdots \longrightarrow V^1_j\overset{f^1}\longrightarrow V^0_j\longrightarrow0$$ is exact. By assumption that $B$ is projective as a graded right $S$-module, we obtain an exact sequence of graded $B$-modules
\begin{equation}\label{eq21}
\cdots\longrightarrow B\ot_S V^i_j\overset{\id\ot_S f^i}\longrightarrow\cdots \longrightarrow B\ot_S V^1_j\overset{\id\ot_S f^1}\longrightarrow B\ot_S V^0_j\longrightarrow0.
\end{equation}
Comparing the generators of each component, one sees that the complex (\ref{eq21}) is isomorphic to the complex $F_jW^\bullet/F_{j-1}W^\bullet$. Therefore, for each $j\ge0$, $F_jW^\bullet/F_{j-1}W^\bullet$ is an exact complex. Now from the classical spectral sequence associated to the filtration of $W^\bullet$ constructed above, we see that $W^\bullet$ is exact. \qed

\begin{prop}\label{prop8} With the same assumptions as in Proposition \ref{prop7}, we have an exact sequence of graded $B$-bimodules:
\begin{equation}\label{eq22}
  \cdots \longrightarrow B\ot_S K_n\ot_S B \overset{d^{-n}}\longrightarrow\cdots\longrightarrow B\ot_S K_1\ot_SB\overset{d^{-1}}\longrightarrow B\ot _SB\overset{\mu}\longrightarrow B\longrightarrow0,
\end{equation}
where the differential $d^{-n}$ is the restriction of the graded $B$-bimodule morphism $B\ot_SM^{\ot_S\,n}\ot_SB\longrightarrow B\ot_SM^{\ot_S\,n-1}\ot_SB$ defined by
\begin{eqnarray*}
 b_1\ot_S x_1\ot_S\cdots\ot_S x_n\ot_S b_2&\mapsto& b_1 x_1\ot_S x_2\ot_S\cdots\ot_S x_n\ot_S b_2\\
&& +(-1)^{n}b_1\ot_S x_1\ot_S\cdots\ot_S x_{n-1}\ot_S x_nb,
\end{eqnarray*}
for $x_1,\dots,x_n\in M,b_1,b_2\in B$. The last morphism $\mu$ is the multiplication of $B$.
\end{prop}
\proof The exactness of the complex (\ref{eq22}) follows from  Lemma \ref{lem5} after applying the functor $S\ot_B-$ on the aforementioned complex. \qed

\section{PBW deformations of Koszul algebras with $pdim(S_B)=2$}\label{sec}
\subsection{Graded deformations}\label{subsec0}
Let $B=B_0\op B_1\op\cdots$ be a locally finite graded algebra. A {\it PBW deformation} of $B$ is a filtered algebra $U$ with an ascending filtration $0=F_{-1}U\subseteq F_0U\subseteq F_1U\subseteq\cdots$ such that the associated graded algebra $gr(U)$ is isomorphic to $B$. Recall from \cite[Remark 1.4]{BG} that a PBW deformation of $B$ is governed by a graded deformation of $B$. Let us recall some notions and terminology of graded deformations.

Let $\k[t]$ be the polynomial algebra with an indeterminate $t$. We view $\k[t]$ as a graded algebra by assuming that the degree of $t$ is 1. Let $B$ be a positively graded algebra. A {\it graded deformation} of $B$ over $\k[t]$ is an associative $\k[t]$-algebra with underlying vector space $B[t]=B\ot \k[t]$ and multiplication defined by
\begin{equation}\label{eq15}
  b_1*b_2=b_1b_2+\mu_1(b_1\ot b_2)t+\mu_2(b_1\ot b_2)t^2+\cdots,
\end{equation}
where $b_1,b_2$ are homogeneous elements in $B$, $b_1b_2$ is the product of $b_1$ and $b_2$ in $B$, and $\mu_i:B\ot B\to B$ is a graded $\k$-linear map of degree $-i$ for all $i\ge1$.

An {\it $i$-th level graded deformation} of $B$ is a graded associative algebra over the graded ring $\k[t]/(t^{i+1})$ with underlying vector space $B[t]/(t^{i+1})$ and multiplication defined by
\begin{equation}\label{eq19}
  b_1*b_2=b_1b_2+\mu_1(b_1\ot b_2)t+\mu_2(b_1\ot b_2)t^2+\cdots+\mu_i(b_1\ot b_2)t^i.
\end{equation}
We say that an $i$-th level graded deformation {\it lifts} to $i+1$-th level graded deformation if there is a $\mu_{i+1}$ such that the multiplication $b_1*b_2=b_1b_2+\mu_1(b_1\ot b_2)t+\cdots+\mu_{i+1}(b_1\ot b_2)t^{i+1}$ defines an associative algebra structure on $B[t]/(t^{i+2})$. If $i$ goes infinite, the obtained graded algebra is called {\it a graded deformation of $B$}, and is usually denoted by $B_t$.

Let $B_t$ be a graded deformation of $B$. Then the associativity of the multiplication of $B_t$ implies the following identities: for $b_1,b_2,b_3\in B$,
\begin{eqnarray}
\label{eq17}  && b_1\mu_1(b_2\ot b_3)-\mu_1(b_1b_2\ot b_3)+\mu_1(b_1\ot b_2b_3)-\mu_1(b_1\ot b_2)b_3=0;\\
\label{eq18}  &\text{for $i\ge2$, }& \sum_{j=1}^{i-1}\mu_j(\mu_{i-j}(b_1\ot b_2)\ot b_3)-\mu_j(b_1\ot \mu_{i-j}(b_2\ot b_3))\\
\nonumber &&=b_1\mu_i(b_2\ot b_3)-\mu_i(b_1b_2\ot b_3)+\mu_i(b_1\ot b_2b_3)-\mu_i(b_1\ot b_2)b_3.
\end{eqnarray}
For $i\ge2$, the linear map $\mathcal{O}^i:B\ot B\ot B\longrightarrow B$, defined by $$\mathcal{O}^i(b_1\ot b_2\ot b_3)=\sum_{j=1}^{i-1}\mu_j(\mu_{i-j}(b_1\ot b_2)\ot b_3)-\mu_j(b_1\ot \mu_{i-j}(b_2\ot b_3)),$$ for $b_1,b_2,b_3\in B$, is called the $(i-1)$-{\it obstruction}.

Let us recall the definition of the bar resolution of $B$. For $n\ge0$, let $Q^{-n}=B\ot B^{\ot n}\ot B$. Define a differential $D^{-n}:Q^{-n}\to Q^{-n+1}$ by $$D^{-n}(b_0\ot b_1\ot\cdots\ot b_n\ot b_{n+1})=\sum_{i=0}^{n}(-1)^{i}b_0\ot\cdots\ot b_ib_{i+1}\ot\cdots b_{n+1},$$ for $b_0,\dots,b_{n+1}\in B$. The complex $(Q^\bullet,D^{\bullet})$ is called {\it the bar resolution} of $B$.
Denote by $CH^n(B,B):=\underline{\Hom}_{B^e}(Q^{-n},B)=\underline{\Hom}(B^{\ot n},B)$ and ${D^n}^*=\underline{\Hom}_{B^e}(D^n,B)$, where $\underline{\Hom}_{B^e}$ stands for the functor of graded $B$-bimodule morphisms, and $\underline{\Hom}$ for the functor of graded morphisms of graded vector spaces. Then $(CH^\bullet(B,B),{D^\bullet}^*)$ is a cochain complex. Recall from \cite{BG}, we have the following description of graded deformations.

\begin{prop}\label{prop6}\cite[Proposition 1.5]{BG} {\rm(i)} A graded linear map $\mu_1:B\ot B\longrightarrow B$ of degree $-1$ determines a first level graded deformation if and only if $\mu_1$ is a 2-cocycle of the complex $(CH^\bullet(B,B),{D^\bullet}^*)$.

{\rm(ii)} Assume that graded linear maps $\mu_k:B\ot B\longrightarrow B$ ($1\leq k\leq i$) define an $i$-th level graded deformation. Then, the $i$-obstruction $\mathcal{O}^i$ is a 3-cocyle of the complex $(CH^\bullet(B,B),{D^\bullet}^*)$.

{\rm(iii)} An $i$-th level graded deformation lifts to an $i+1$-th graded deformation if and only if there is a graded linear map $\mu_{i+1}:B\ot B\to B$ of degree $-i-1$ such that $\mathcal{O}^{i+1}={D^{-2}}^*(\mu_{i+1})$.
\end{prop}

\subsection{Graded projective resolutions of Koszul algebras with pdim$S_B$=2} \label{subsec2}
Now assume $B$ is a Koszul algebra. Like before, we write $S$ for $B_0$ and $M$ for $B_1$. Let $B\cong T_S(M)/(R)$ for some sub-$S$-bimodule $R\subseteq M\ot_S M$.

From now on till the end of this section, we assume that the projective dimension of $S_B$ is two.

Since $B$ is Koszul and $S_B$ is of graded projective dimension two, we have $M\ot_S R\cap R\ot_S M=0$ by Proposition \ref{prop7}. It follows from Proposition \ref{prop8} that the following sequence of graded $B$-bimodules is exact:
\begin{equation}\label{eq24}
   0\longrightarrow B\ot_S R\ot_S B \overset{d^{-2}}\longrightarrow B\ot_S M\ot_SB\overset{d^{-1}}\longrightarrow B\ot _SB\overset{\mu}\longrightarrow B\longrightarrow0.
\end{equation}
This sequence is not a projective resolution of the graded $B$-bimodule ${}_BB_B$ in general. However, we may construct a projective resolution from this sequence.

We next construct $B$-bimodule projective resolutions of the first three nonzero terms in the above sequence. We need the bar resolution of the $S$-module $R$, $M$ and $S$ respectively. Recall that the bar resolution of a left $S$-module $N$ is defined as follows:
\begin{equation}\label{eq25}
   \cdots\longrightarrow S\ot S\ot S\ot N  \overset{d_N^{-2}}\longrightarrow S\ot S\ot N\overset{d_N^{-1}}\longrightarrow S\ot N\overset{\mu_N}\longrightarrow N\longrightarrow0,
\end{equation}
where the differential is given by
\begin{eqnarray*}
  d_N^{-i}(s_0\ot s_1\ot \dots\ot s_i\ot x)&=&\sum_{j=0}^{i-1}(-1)^ks_0\ot\cdots\ot s_ks_{k+1}\ot\cdots\ot s_i\ot x\\
  &&+(-1)^is_0\ot\cdots\ot s_{i-1}\ot s_ix
\end{eqnarray*}
for all $i\ge0$ and $s_0,\dots,s_i\in S$ and $x\in N$. We denote by $\mathcal{B}(N)$  the bar resolution of $N$. If $N$ is an $S$-bimodule, then the differential $d_N$ of the bar resolution is compatible with the right $S$-module structure of each component.

Assume further that $N$ is an $S$-bimodule which is both left and right projective. Applying the functors $B\ot_S-$ and $-\ot_S B$ simultaneously on the bar resolution $\mathcal{B}(N)$, we obtain a complex
\begin{equation}\label{eq26}
   \cdots\longrightarrow B\ot S\ot S\ot N\ot_S B  \overset{\partial^{-2}}\longrightarrow B\ot S\ot N\ot_SB\overset{\partial^{-1}}\longrightarrow B\ot N\ot_S B\overset{\partial^0}\longrightarrow B\ot_SN\ot_SB\longrightarrow0,
\end{equation} where the differential $\partial^{-i}$ is induced by the differential of the bar resolution of $N$.
Since $N$ is both left and right projective, the sequence above is exact. Moreover, each component $B\ot S^{\ot i}\ot N\ot_SB$ is projective as a $B$-bimodule. Hence the sequence (\ref{eq26}) provides a projective resolution of the $B$-bimodule $B\ot_SN\ot_SB$.

Now taking $N$ to be $R$, $M$ and $S$ respectively in the sequence (\ref{eq26}), we obtain projective resolutions for the first three non-zero terms of the sequence (\ref{eq24}) respectively. Comparing the projective resolutions of each components, we obtain a commutative diagram (ignoring the diagonal arrows):
{\small
$$\xymatrix{
 0&0&0\\
 B\ot_S R\ot_S B\ar[u]\ar[r]^{d^{-2}}&B\ot_S M\ot_SB\ar[u]\ar[r]^{d^{-1}}&B\ot_SB\ar[u]\\
  B\ot R\ot_S B\ar[drr]^{h^0\qquad\qquad} \ar[u]_{\partial^{0}} \ar[r]^{\vartheta^{-2,0}} & B\ot M\ot_S B \ar[u]_{\partial^0} \ar[r]^{\vartheta^{-1,0}} & B\ot B\ar[u]_{\partial^0}\\
  B\ot  S\ot R\ot_S B\ar[drr]^{h^{-1}\qquad\qquad} \ar[u]_{\partial^{-1}} \ar[r]^{\vartheta^{-2,-1}} & B\ot S\ot M\ot_S B \ar[u]_{\partial^{-1}} \ar[r]^{\vartheta^{-1,-1}} & B\ot S\ot B \ar[u]_{\partial^{-1}} \\
    B\ot S\ot S\ot R\ot_S B\ar[drr]^{h^{-2}\qquad\qquad} \ar[u]_{\partial^{-2}} \ar[r]^{\vartheta^{-2,-2}} & B\ot S\ot S\ot M\ot_S B \ar[u]_{\partial^{-2}} \ar[r]^{\vartheta^{-1,-2}} & B\ot S\ot S\ot B \ar[u]_{\partial^{-2}}\\
    \vdots\ar[u]&\vdots\ar[u]&\vdots\ar[u] }$$
$$\text{Diagram 1}$$
}
where the morphisms $\vartheta^{-2,-i}$ ($i\ge0$) are induced by $d^{-2}$, and the morphisms $\vartheta^{-1,-j}$ ($j\ge0$) induced by $d^{-1}$. In general, the composition $\vartheta^{-1,-i}\vartheta^{-2,-i}$ is not zero. Since $d^{-1}d^{-2}=0$, the cochain complex morphism $\vartheta^{-1,\bullet}\vartheta^{-2,\bullet}$ is homotopic to zero. Hence there are graded $B$-bimodules morphisms $$h^{-i}:B\ot S^{\ot i}\ot R\ot_S B\longrightarrow B\ot S^{\ot i+1}\ot  B$$ for all $i\ge 0$ such that $\vartheta^{-1,-i}\vartheta^{-2,-i}=\partial^{-i-1}h^{-i}+h^{-i+1}\partial^{-i}$ for all $i\ge1$ and $\vartheta^{-1,0}\vartheta^{-2,0}=\partial^{-1}h^{0}$.
Now let $$P^{-n}= B\ot S^{\ot n-2}\ot R\ot_SB\bigoplus B\ot S^{\ot n-1}\ot M\ot_SB\bigoplus B\ot S^{\ot n}\ot B \text{\ for $n\ge2$},$$
$$P^{-1}=B\ot M\ot_S B\bigoplus B\ot S\ot B,\text{ and } P^0=B\ot B.$$
Define morphisms $\delta^{-n}:P^{-n}\to P^{-n+1}$ by
\begin{eqnarray*}
  \delta^{-n}(x)&=&\partial^{-n+2}(x)+\vartheta^{-2,-n+2}(x)-h^{-n+2}(x)\text{ for all $n\ge3$ and $x\in B\ot S^{\ot n-2}\ot R\ot_SB$}, \\ \delta^{-n}(y)&=&-\partial^{-n+1}(y)+\vartheta^{-1,-n+1}(y)\text{ for all $n\ge2$ and $y\in B\ot S^{\ot n-1}\ot M\ot_SB$}, \\ \delta^{-n}(z)&=&-\partial^{-n}(z)\text{ for all $n\ge1$ and $z\in B\ot S^{\ot n}\ot B$},\\
   \delta^{-2}(x)&=&\vartheta^{-2,0}(x)-h^0(x)\text{ for all $x\in B\ot R\ot_SB$},\\
   \delta^{-1}(y)&=&\vartheta^{-1,0}(y)\text{ for $y\in B\ot M\ot_SB$}.
\end{eqnarray*}
Now it is easy to check $\delta^{-n}\delta^{-n-1}=0$ for all $n\ge1$. Hence $(P^\bullet,\delta^\bullet)$ forms a complex. Take the subcomplex $F_0P^\bullet$ to be the 0th column of Diagram 1 above, $F_{-1}P^\bullet$ to be the subcomplex consists of 0th and $(-1)$th column, $F_{-i}P^\bullet=P^\bullet$ for all $i\ge2$, and $F_jP^\bullet=0$ for all $j\ge1$. Then $\cdots\supseteq F_{-i}P^\bullet\supseteq F_{-i+1}P^\bullet\supseteq\cdots$ is an exhaustive filtration of the complex $P^\bullet$. The first level of the spectral sequence associated to this filtration collapses, and we see that $P^\bullet$ is exact except at the 0th position, and the cohomology of 0th position is $H^0P^\bullet= B$. Therefore we obtain the following result. Note that if $S$ is a group algebra (or Hopf algebra) and $A$ is a module algebra over $S$, then the projective resolution $P^\bullet$ of the smash product $A\#S$ is equivalent to the resolution obtained in \cite{SW} (or \cite{WW}, also \cite{Ne}) by choosing suitable $\vartheta$ and $h$.

\begin{prop}\label{prop9} With the notions as above. $(P^\bullet,\delta^\bullet)$ is a graded projective resolution of the $B$-bimodule $B$.
\end{prop}

\subsection{Comparing the projective resolution and the bar resolution of $B$}\label{subsec1}

Note that each component of $P^\bullet$ in Proposition \ref{prop9} is generated in degrees 0, or 1, or 2. Hence the third Hochschild cohomology $HH^3(B,B)$ is concentrated in degrees not less than $-2$. If we are given an $i$th level of graded deformation of $B$, then the $i$-obstruction $\mathcal{O}^i$ must be a coboundary and hence the graded deformations automatically lifts to the next level. To study the graded deformations of $B$ more closely, we need to compare the projective resolution $P^\bullet$ and the bar resolution $Q^\bullet$ of $B$. Let us first construct the morphisms $\vartheta^{-2,0}$, $\vartheta^{-1,0}$ and $h^0$ in Diagram 1 more concretely. We need more notations.

Let $X$ be a right projective $S$-module.
Let $\mu_{{}_X}:X\ot S\to X$ be the right $S$-action on $X$ which is obviously a right $S$-module morphism. Since $X$ is right projective, the right $S$-module morphism $\mu_{{}_X}$ splits. Hence there is a right $S$-module morphism $\rho_{{}_X}:X\to X\ot S$ such that $\mu_{{}_X}\rho_{{}_X}=\id$. To simplify the notions, we use Sweedler's notion \cite{S} to denote the actions of $\rho_{{}_X}$ on the elements in $X$, that is, for $x\in X$, we write $\rho_{{}_X}(x)=x^{(0)}\ot x^{(1)}$.
\begin{lem}\label{lem6} With the notions above, we have the following properties:

{\rm(i)} for $x\in X$, $x^{(0)}x^{(1)}=x$;

{\rm(ii)} for $x\in X$ and $s\in S$, $(xs)^{(0)}\ot (xs)^{(1)}=x^{(0)}\ot x^{(1)}s$;

{\rm(iii)} for $x\in X$, $\left(x^{(0)}\right)^{(0)}\ot \left(x^{(0)}\right)^{(1)}x^{(1)}=x^{(0)}\ot x^{(1)}$.
\end{lem}
\proof The property (i) holds because $\mu_{{}_X}\rho_{{}_X}=\id$. The property (ii) holds because $\rho_{{}_X}$ is a right $S$-module morphism. The property (iii) follows from (i) and (ii).\qed

\begin{rem} \label{rem2} {\rm Since $B$ is both left and right projective, for each $n\ge1$, we fix an $S$-module morphism $\rho_{{}_{B_n}}:B_n\to B_n\ot S$ such that $\mu_{{}_{B_n}}\rho_{{}_{B_n}}=\id$, and we always use the notion introduced above to denote the action of $\rho_{{}_{B_n}}$ on elements in the following narratives.}
\end{rem}

Now we are going to construct the morphisms $\vartheta^{-2,0}$ and $\vartheta^{-1,0}$ as appeared in the Diagram 1. At first, let us consider the following exact sequence of $S$-bimodules $$0\longrightarrow R\overset{\iota}\longrightarrow M\ot_S M\overset{\mu}\longrightarrow B_2\longrightarrow0,$$ where $\iota$ is the inclusion map and $\mu$ is the multiplication of $B$. Since $B$ is right projective, the sequence above splits. Hence there is a right $S$-module morphism $\xi: M\ot_SM\to R$ such that $\xi\iota=\id$. Note that the natural projection map $\pi:M\ot M\to M\ot_SM$ is also a right $S$-module morphism. Hence $\pi$  also splits. Indeed, if we let $\zeta$ be the composition
\begin{equation}\label{eq30}
\zeta:M\ot_SM\overset{\rho_M\ot_S\id}\longrightarrow M\ot S\ot_SM\cong M\ot M,
\end{equation}
then we have $\pi\zeta=\id$.
Now let $\alpha$ be the composition of the right $S$-module morphisms:
$$\xymatrix{\alpha: R\ar[r]^{\iota} &M\ot_SM\ar[r]^{\zeta} &M\ot M\ar[r]^{\id\ot \rho_M\quad}& M\ot M\ot S.  }$$
For $r=\sum_{i=1}^nx_i\ot_Sy_i\in M\ot_SM$, we see $$\alpha(r)=\sum_{i=1}^nx_i^{(0)}\ot \left(x_i^{(0)}y_i\right)^{(0)}\ot \left(x_i^{(0)}y_i\right)^{(1)}.$$
With the above notions, we define $B$-bimodule morphisms $\vartheta^{-2,0}: B\ot R\ot_SB\longrightarrow B\ot M\ot_SB$ by
\begin{equation}\label{eq27}
  \vartheta^{-2,0}(1\ot r\ot_S 1)=\sum_{i=1}^nx_i^{(0)}\ot x_i^{(1)}y_i\ot_S 1+\sum_{i=1}^n 1\ot x_i\ot_S y_i,
\end{equation}
and $\vartheta^{-1,0}:B\ot M\ot_SB\longrightarrow B\ot B$ by
\begin{equation}\label{eq28}
  \vartheta^{-1,0}(1\ot m\ot_S 1)=m^{(0)}\ot m^{(1)}-1\ot m.
\end{equation}
The morphisms are well defined since $\rho_M$ is a right $S$-module morphism. It is straightforward to check that $\partial^0\vartheta^{-2,0}=d^{-2}\partial^0$ and $\partial^0\vartheta^{-1,0}=d^{-1}\partial^0$.

We next construct the morphism $h^0$ as appeared in Diagram 1. Let $\overline{h}$ be the composition of the following right $S$-module morphisms: $$\overline{h}:R\overset{\alpha}\longrightarrow M\ot M\ot S\overset{f}\longrightarrow B\ot M\ot M\ot B\overset{D^{-1}}\longrightarrow B\ot (M\op B_2)\ot B\overset{g}\longrightarrow B\ot S\ot B, $$ where $D^{-1}$ is the differential of the bar resolution of $B$, $f$ is defined by $f(m_1\ot m_2\ot s)=1\ot m_1\ot m_2\ot s$ for $m_1,m_2\in M$ and $s\in S$, and $g$ is the $B$-bimodule morphism defined by $g(1\ot x\ot 1)=x^{(0)}\ot x^{(1)}\ot 1$ for all $x\in M$ or $x\in B_2$. Let $$h^0:B\ot R\ot_SB\longrightarrow B\ot S\ot B$$ be the $B$-bimodule morphism determined by $-\overline{h}$. For $r=\sum_{i=1}^nx_i\ot_Sy_i\in R$, we may check that $h^0$ acts on element $1\ot r\ot_S 1$  as follows:
{\small\begin{eqnarray*}
  &&h^0(1\ot r\ot_S 1)\\
  &=&-g\sum_{i=1}^n\left[x_i^{(0)}\ot \left( x_i^{(1)}y_i\right)^{(0)}\ot \left( x_i^{(1)}y_i\right)^{(1)}
  -1\ot x_i^{(0)} \left( x_i^{(1)}y_i\right)^{(0)}\ot \left( x_i^{(1)}y_i\right)^{(1)}+ 1\ot x_i^{(0)}\ot x_i^{(1)}y_i\right]\\
  &=& -\sum_{i=1}^n x_i^{(0)}\left[ \left( x_i^{(1)}y_i\right)^{(0)}\right]^{(0)}\ot\left[ \left( x_i^{(1)}y_i\right)^{(0)}\right]^{(1)} \ot \left( x_i^{(1)}y_i\right)^{(1)}\hspace{25mm} \text{(I)}\\
  &&+\sum_{i=1}^n\left[x_i^{(0)} \left( x_i^{(1)}y_i\right)^{(0)}\right]^{(0)}\ot\left[x_i^{(0)} \left( x_i^{(1)}y_i\right)^{(0)}\right]^{(1)} \ot \left( x_i^{(1)}y_i\right)^{(1)}\hspace{18mm} \text{(II)}\\
  &&-\sum_{i=1}^n\left(x_i^{(0)}\right)^{(0)}\ot \left(x_i^{(0)}\right)^{(1)}\ot x_i^{(1)}y_i. \hspace{55mm} \text{(III)}
\end{eqnarray*}}
We next check that $\partial^{-1}h^0=\vartheta^{-1,0}\vartheta^{-2,0}$. We have the following computations in which we freely use the properties in Lemma \ref{lem6} and the fact that $\sum_{i=1}^n x_i y_i=0$:
\begin{eqnarray*}
 &&\vartheta^{-1,0}\vartheta^{-2,0}(1\ot r\ot_S 1)\\
 &=& \vartheta^{-1,0}\left(\sum_{i=1}^nx_i^{(0)}\ot x_i^{(1)}y_i\ot_S 1+\sum_{i=1}^n 1\ot x_i\ot_S y_i\right)\\
 &=&\sum_{i=1}^nx_i^{(0)}\left[x_i^{(1)}y_i\right]^{(0)}\ot \left[x_i^{(1)}y_i\right]^{(1)}-\sum_{i=1}^nx_i^{(0)}\ot x_i^{(1)}y_i+\sum_{i=1}^nx_i^{(0)}\ot x_i^{(1)}y_i\\
 &=&\sum_{i=1}^nx_i^{(0)}\left[x_i^{(1)}y_i\right]^{(0)}\ot \left[x_i^{(1)}y_i\right]^{(1)}.
\end{eqnarray*}
Let us compute $\partial^{-1}h^0(1\ot r\ot_S1)$. As we see above that $h^0(1\ot r\ot_S1)$ has three terms (I), (II) and (III). Applying $\partial^{-1}$ on these terms individually, we obtain
\begin{eqnarray*}
&& \partial^{-1}\left(\sum_{i=1}^n x_i^{(0)}\left[ \left( x_i^{(1)}y_i\right)^{(0)}\right]^{(0)}\ot\left[ \left( x_i^{(1)}y_i\right)^{(0)}\right]^{(1)} \ot \left( x_i^{(1)}y_i\right)^{(1)}\right)\\
&=& \sum_{i=1}^n x_i^{(0)}\left[ \left( x_i^{(1)}y_i\right)^{(0)}\right]^{(0)}\left[ \left( x_i^{(1)}y_i\right)^{(0)}\right]^{(1)} \ot \left( x_i^{(1)}y_i\right)^{(1)}\\
&&-\sum_{i=1}^n x_i^{(0)}\left(\left[ \left( x_i^{(1)}y_i\right)^{(0)}\right]^{(0)}\ot\left[ \left( x_i^{(1)}y_i\right)^{(0)}\right]^{(1)}\left( x_i^{(1)}y_i\right)^{(1)}\right)\\
&=& \sum_{i=1}^n x_i^{(0)} \left( x_i^{(1)}y_i\right)^{(0)} \ot \left( x_i^{(1)}y_i\right)^{(1)}
-x_i^{(0)} \left( x_i^{(1)}y_i\right)^{(0)} \ot \left( x_i^{(1)}y_i\right)^{(1)}\\
&=&0,
\end{eqnarray*}
\begin{eqnarray*}
 && \partial^{-1}\left(\sum_{i=1}^n \left(x_i^{(0)}\right)^{(0)}\ot \left(x_i^{(0)}\right)^{(1)}\ot x_i^{(1)}y_i\right)\\
 &=&\sum_{i=1}^n\left(x_i^{(0)}\right)^{(0)} \left(x_i^{(0)}\right)^{(1)}\ot x_i^{(1)}y_i-\left(x_i^{(0)}\right)^{(0)}\ot \left(x_i^{(0)}\right)^{(1)} x_i^{(1)}y_i\\
 &=&\sum_{i=1}^n x_i^{(0)}\ot x_i^{(1)}y_i-\left[\left(x_i^{(0)}\right)^{(0)}\ot\left(x_i^{(0)}\right)^{(1)} x_i^{(1)}\right]y_i\\
 &=&\sum_{i=1}^n x_i^{(0)}\ot x_i^{(1)}y_i-x_i^{(0)}\ot x_i^{(1)}y_i\\
 &=&0,\end{eqnarray*} and
\begin{eqnarray*}
 &&\partial^{-1}\left(\sum_{i=1}^n\left[x_i^{(0)} \left( x_i^{(1)}y_i\right)^{(0)}\right]^{(0)}\ot\left[x_i^{(0)} \left( x_i^{(1)}y_i\right)^{(0)}\right]^{(1)} \ot \left( x_i^{(1)}y_i\right)^{(1)}\right)\\
 &=&\sum_{i=1}^n \left[x_i^{(0)} \left( x_i^{(1)}y_i\right)^{(0)}\right] \ot \left( x_i^{(1)}y_i\right)^{(1)}\\
 &&-\sum_{i=1}^n\left[x_i^{(0)} \left( x_i^{(1)}y_i\right)^{(0)}\right]^{(0)}\ot\left[x_i^{(0)} \left( x_i^{(1)}y_i\right)^{(0)}\right]^{(1)} \left( x_i^{(1)}y_i\right)^{(1)}\\
 &=&\sum_{i=1}^n x_i^{(0)} \left( x_i^{(1)}y_i\right)^{(0)} \ot \left( x_i^{(1)}y_i\right)^{(1)}.
\end{eqnarray*}
Hence $\partial^{-1}h^0=\vartheta^{-1,0}\vartheta^{-2,0}$ as desired.

\subsection{Construction of a specific homotopic equivalence from $Q^\bullet$ to $P^\bullet$}\label{subsec3}

Next we construct a specific homotopic equivalence $\varphi^\bullet:Q^\bullet\to P^\bullet$. Since $Q^0=P^0=B\ot B$, we define $\varphi^0:Q^0\to P^0$ to be the identity map. Since $Q^{-1}=B\ot B\ot B=B\ot M\ot B\bigoplus B\ot \left(\op_{n\ge2}B_n\right)\ot B$. We define $\varphi^{-1}:Q^{-1}\to P^{-1}$ as follows. For $m\in M$, we let
\begin{equation}\label{eq29}
  \varphi^{-1}(1\ot m\ot 1)=1\ot m\ot_S 1+m^{(0)}\ot m^{(1)}\ot 1.
\end{equation}
Note that the element $1\ot m\ot_S 1$ is in $B\ot M\ot_S B$ and the element $m^{(0)}\ot m^{(1)}\ot 1$ is in $B\ot S\ot B$.

For the restriction of $\varphi^{-1}$ to $B\ot \left(\op_{n\ge2}B_n)\right)\ot B$, we define $\varphi^{-1}|_{B\ot \left(\op_{n\ge2}B_n\right)\ot B}$ to be any graded $B$-bimodule morphism satisfying the identity $\varphi^0 D^{-1}(x)=\delta^{-1}\varphi^{-1}(x)$ for all $x\in B\ot \left(\op_{n\ge2}B_n\right)\ot B$. Now for all $m\in M$, we have
\begin{eqnarray*}
 \delta^{-1}\varphi^{-1}(1\ot m\ot 1)&=&\delta^{-1}(1\ot m\ot_S 1+m^{(0)}\ot m^{(1)}\ot 1)\\
 &=&m^{(0)}\ot m^{(1)}-1\ot m+m\ot 1-m^{(0)}\ot m^{(1)}\\
 &=&m\ot 1-1\ot m\\
 &=&\varphi^0 D^{-1}(1\ot m\ot 1).
\end{eqnarray*}
Hence $\delta^{-1}\varphi^{-1}=\varphi^0 D^{-1}$ holds on $B\ot M\ot B$, and therefore it holds on $Q^{-1}$.

Now we construct $\varphi^{-2}: Q^{-2}\to P^{-2}$. Let $\beta$ be the composition of the following morphisms
$$\beta:R\overset{\iota}\longrightarrow M\ot_SM\overset{\zeta}\longrightarrow M\ot M,$$ where $\iota$ is the inclusion map and $\zeta$ is defined by (\ref{eq30}) in Subsection \ref{subsec1}. Let $\gamma$ be the composition of the following morphism
$$\gamma:M\ot M\overset{\pi}\longrightarrow M\ot_SM\overset{\xi}\longrightarrow R,$$ where the morphisms $\pi$ and $\xi$ are defined in Subsection \ref{subsec1}.
Since $\pi\zeta=\id$ and $\xi\iota=\id$ (cf. Subsection \ref{subsec1}), we have $\gamma\beta=\id$. Hence we have a decomposition $$M\ot M=\text{im}(\beta)\op \ker(\gamma).$$
Therefore,
 $$B\ot M\ot M\ot B=B\ot \text{im}(\beta)\ot B\bigoplus B\ot\ker(\gamma)\ot B.$$
Finaly, we have the following decomposition of graded $B$-bimodules $$Q^{-2}=B\ot M\ot M\ot B\bigoplus K=B\ot \text{im}(\beta)\ot B\bigoplus B\ot\ker(\gamma)\ot B\bigoplus K,$$ where $K$ is a free complement of $B\ot M\ot M\ot B$ in $Q^{-2}$. For the restriction of $\varphi^{-2}$ to the direct summands $B\ot\ker(\gamma)\ot B\bigoplus K$, we may choose any graded $B$-bimodule morphism such that the identity $\delta^{-2}\circ\varphi^{-2}_{B\ot\ker(\gamma)\ot B\bigoplus K}=\varphi^{-1}\circ D^{-2}$ holds on $B\ot\ker(\gamma)\ot B\bigoplus K$.

We define the restriction of $\varphi^{-2}$ to $B\ot \text{im}(\beta)\ot B$ as follows. Let $\overline{\varphi}:B\ot \text{im}(\beta)\ot B\longrightarrow B\ot R\ot_S B$ be the graded $B$-module morphism determined by the right $S$-module morphism $\gamma:\text{im}(\beta)\longrightarrow R$. Let us check the difference between the composition maps $\delta^{-2}\overline{\varphi}$ and $\varphi^{-1}D^{-2}$. For $r=\sum_{i=1}^nx_i\ot_Sy_i\in R$, we have $\beta(r)=\sum_{i=1}^nx_i^{(0)}\ot x_i^{(1)}y_i$. Hence
\begin{eqnarray*}
 &&\left[\delta^{-2}\overline{\varphi}-\varphi^{-1}D^{-2}\right](1\ot \beta(r)\ot 1)\\
 &=&\delta^{-2}(1\ot r\ot_S 1)-\varphi^{-1}\left(\sum_{i=1}^nx_i^{(0)}\ot x_i^{(1)}y_i\ot 1+\sum_{i=1}^n1\ot x_i^{(0)}\ot x_i^{(1)}y_i\right)\\
 &=&-h^0(1\ot r\ot_S1)-\left[\sum_{i=1}^nx_i^{(0)}\ot x_i^{(1)}y_i\ot 1+\sum_{i=1}^n\left(x_i^{(0)}\right)^{(0)}\ot \left(x_i^{(0)}\right)^{(1)}\ot x_i^{(1)}y_i\right].
\end{eqnarray*}
Note that both $h^0(1\ot r\ot_S1)$ and $\sum_{i=1}^nx_i^{(0)}\ot x_i^{(1)}y_i\ot 1+\sum_{i=1}^n\left(x_i^{(0)}\right)^{(0)}\ot \left(x_i^{(0)}\right)^{(1)}\ot x_i^{(1)}y_i$ are in $B\ot S\ot B$. From the computations in Subsection \ref{subsec1}, we have $$\partial^{-1}\left(h^0(1\ot r\ot_S1)\right)=\sum_{i=1}^nx_i^{(0)}\left(x_i^{(1)}y_i\right)^{(0)}\ot \left(x_i^{(1)}y_i\right)^{(1)}.$$ On the other hand,
\begin{eqnarray*}
 &&\partial^{-1}\left[\sum_{i=1}^nx_i^{(0)}\ot x_i^{(1)}y_i\ot 1+\sum_{i=1}^n\left(x_i^{(0)}\right)^{(0)}\ot \left(x_i^{(0)}\right)^{(1)}\ot x_i^{(1)}y_i\right]\\
 &=& -\sum_{i=1}^nx_i^{(0)}\left(x_i^{(1)}y_i\right)^{(0)}\ot \left(x_i^{(1)}y_i\right)^{(1)}+\sum_{i=1}^nx_i^{(0)}\ot x_i^{(1)}y_i-\sum_{i=1}^n\left(x_i^{(0)}\right)^{(0)}\ot \left(x_i^{(0)}\right)^{(1)} x_i^{(1)}y_i\\
 &=&-\sum_{i=1}^nx_i^{(0)}\left(x_i^{(1)}y_i\right)^{(0)}\ot \left(x_i^{(1)}y_i\right)^{(1)}+\sum_{i=1}^nx_i^{(0)}\ot x_i^{(1)}y_i-\sum_{i=1}^nx_i^{(0)}\ot x_i^{(1)}y_i\\
 &=&-\sum_{i=1}^nx_i^{(0)}\left(x_i^{(1)}y_i\right)^{(0)}\ot \left(x_i^{(1)}y_i\right)^{(1)}.
\end{eqnarray*}
Therefore $\partial^{-1}\left[\delta^{-2}\overline{\varphi}-\varphi^{-1}D^{-2}\right](1\ot \beta(r)\ot 1)=0$, which implies that the image of the morphism $\delta^{-2}\overline{\varphi}-\varphi^{-1}D^{-2}:B\ot \text{im}(\beta)\ot B\longrightarrow B\ot S\ot B$ lies in $\ker(\partial^{-1})=\text{im}(\partial^{-2})$. Hence there is a graded $B$-bimodule morphism $$\eta:B\ot \text{im}(\beta)\ot B\longrightarrow B\ot S\ot S\ot B$$ such that $\partial^{-2}\eta=\delta^{-2}\overline{\varphi}-\varphi^{-1}D^{-2}$. Now let $$\varphi^{-2}|_{B\ot \text{im}(\beta)\ot B}=\overline{\varphi}-\eta.$$
Then on the direct summand $B\ot \text{im}(\beta)\ot B$, we have
\begin{eqnarray*}
 \delta^{-2}\varphi^{-2}&=&\delta^{-2}\overline{\varphi}-\delta^{-2}\eta\\
 &=&\delta^{-2}\overline{\varphi}-\partial^{-2}\eta\\
 &=&\delta^{-2}D^{-2}.
\end{eqnarray*}
As a consequence, $\delta^{-2}\varphi^{-2}=\delta^{-2}D^{-2}$ holds on $Q^{-2}$.

\subsection{PBW deformations of $B$}

Let $B=T_S(M)/(R)$ be a Koszul algebra such that pdim$(S_B)=2$. Let $\theta:R\to S$ be a linear map. Let $U=T_S(M)/(r-\theta(r):r\in R)$ be the quotient algebra of $T_S(M)$ divided by the ideal generated by all the elements $r-\theta(r)$ ($r\in R$). Let $p:T_S(M)\to U$ be the projection map. For $i\ge0$, set $T^{\leq i}_S(M)=S\op M\op\cdots\op M^{\ot_S i}$, and $F_iU=p(T^{\leq i}_S(M))$. Then $U$ is a filtered algebra with the ascending filtration $0\subseteq F_0U\subseteq \cdots \subseteq F_iU\subseteq\cdots$. We have the following main result of this section.

\begin{thm}\label{thm3} $U$ is a PBW deformation of $B$ if and only if $\theta:R\to S$ is an $S$-bimodule morphism.
\end{thm}
\proof Assume $U$ is a PBW deformation of $B$. Let $I=(r-\theta(r):r\in R)$. Then $I\cap T^{\leq1}_S(M)=0$. For all $r\in R$ and $s\in S$, both elements $sr-\theta(sr)$ and $sr-s\theta(r)$ are in $I$. Hence $\theta(sr)-s\theta(r)$ is in $I$. On the other hand $\theta(sr)-s\theta(r)\in S\subseteq T^{\leq1}_S(M)$, which implies $\theta(sr)-s\theta(r)=0$. Therefore $\theta$ is a left $S$-module morphism. Similarly, we see that $\theta$ is a right $S$-module morphism.

Now assume $\theta$ is an $S$-bimodule morphism. Consider the graded projective resolution $P^\bullet$ of ${}_BB_B$ constructed in Subsection \ref{subsec2} with the specific differentials constructed in Subsection \ref{subsec1}. Let $C^n(B,B)=\underline{\Hom}_{B^e}(P^{-n},B)$ and ${\delta^{-n}}^*=\underline{\Hom}_{B^e}(\delta^{-n},B)$. Then $(C^\bullet(B,B),{\delta^\bullet}^*)$ is a cochain complex. Let $\varphi^\bullet:Q^\bullet\to P^\bullet$ be the homotopic equivalence constructed in Subsection \ref{subsec3}. Then ${\varphi^\bullet}^*:(C^\bullet(B,B),{\delta^\bullet}^*)\longrightarrow(CH^\bullet(B,B),{D^\bullet}^*)$ is a homotopic equivalence of the cochain complex. Note that $P^{-2}=B\ot R\ot_SB\bigoplus B\ot S\ot M\ot_S B\bigoplus B\ot S\ot S \ot B$. Define a graded $B$-bimodule morphism $\widehat{\theta}:P^{-2}\longrightarrow B$ of degree $-2$ by $\widehat{\theta}(a\ot r\ot_Sb)=a\theta(r)b$ for all $a,b\in B$ and $r\in R$, so that $\widehat{\theta}(x)=0$ for all $x\in B\ot S\ot M\ot_S\ot B$ or $x\in B\ot S\ot S \ot B$. The morphism $\widehat{\theta}$ is well defined since $\theta$ is a right $S$-module morphism. Since $\theta$ is also a left $S$-module morphism, we see that $\widehat{\theta}$ is a 2-cocycle in the cochain complex $C^\bullet(B,B)$. Now let $\mu_2={\varphi^{-2}}^*(\widehat{\theta})\in CH^2(B,B)$. Since $\varphi^\bullet$ is a homotopic equivalence  respecting the gradings of the modules, we see that $\mu_2$ is a graded $2$-cocycle of degree $-2$. Let $\mu_1=0$. Then $\mu_1$ and $\mu_2$ define a second level of graded deformation of $B$. Note that the obstruction $\mathcal{O}^i$ for a graded deformation of $B$ is of degree $-i$. Since the graded projective $B$-bimodule $P^{-n}$ is generated in degrees 0, 1 and 2, the $n$th cohomology of $C^\bullet(B,B)$ is concentrated in degree not less than $-2$ for all $n\ge0$, which in turn implies that the $n$th Hochschild cohomology $HH^n(B,B)$ is concentrated in degrees not less than $-2$. Hence $\mathcal{O}^i$ must be a coboundary for all $i\ge3$. Therefore, there are graded maps $\mu_i:B\ot B\to B$ of degree $-i$ such that $\mu_1,\mu_2,\mu_3,\dots,$ define a graded deformation of $B$, that is, $B_t=B[t]$ with the multiplication defined by $b_1*b_2=b_1b_2+\mu_1(b_1\ot b_2)t+\mu_2(b_1\ot b_2)t^2+\cdots$ is a graded algebra. Specializing at $t=1$, we obtain a PBW deformation $\overline{U}$ of $B$. We show that $\overline{U}$ is isomorphic to $U$.

Let us inspect the action of $\mu_2$ on $B\ot B$ more closely. Since $\mu_2$ is of degree $-2$, $\mu_2(s\ot m)=\mu_2(m\ot s)=0$ for all $s\in S$ and $m\in M$. Then, from the associativity of $B_t$, we see that $\mu_2(m_1s\ot m_2)=\mu_2(m_1\ot sm_2)$. Hence $\mu_2$ facts through the projection $\pi:M\ot M\to M\ot_SM$, that is, there is a map $\overline{\mu}_2:M\ot_SM\to B$ such that $\mu_2|_{M\ot M}=\overline{\mu}_2\pi$. By definition, $\mu_2={\varphi^{-2}}^*(\widehat{\theta})$. From the definition of $\varphi^{-2}$ in Subsection \ref{subsec3}, we have $\overline{\mu}_2|_R=\theta$.

Now $\overline{U}$ is clearly generated by $S\op M$. Hence we have an algebra epimorphism $\overline{p}:T_S(M)\to \overline{U}$. For $r\in R$, the arguments above show that $\overline{p}(r-\theta(r))=0$. Hence $\overline{p}$ induces an algebra epimorphism $q:U\to \overline{U}$. Let $F_i\overline{U}=\overline{p}(T_S^{\leq i}(M))$ for all $i\ge0$. Clearly, $F_0U=F_0\overline{U}=S$. Since $U=T_S(M)/(r-\theta(r):r\in R)$, we see that $\dim F_i(U)/F_{i-1}U\leq \dim B_i$ for all $i\ge1$. We show inductively that the restriction $q:F_iU\to F_i\overline{U}$ is bijective for all $i\ge0$. Suppose that we already have that the restriction $q:F_nU\to F_n\overline{U}$ is bijective. Since $q$ is an epimorphism, $\dim F_{n+1}U/F_nU\ge \dim F_{n+1}\overline{U}/F_{n+1}\overline{U}$. Since $\overline{U}$ is a PBW deformation of $B$ we have $F_{n+1}\overline{U}/F_{n+1}\overline{U}\cong B_{n+1}$ as vector spaces. Hence we have the following inequality $\dim B_{n+1}\ge \dim F_{n+1}U/F_nU\ge \dim B_{n+1}$, which implies $\dim F_{n+1}U/F_nU=\dim B_{n+1}$. Hence $U\cong \overline{U}$.\qed

\begin{exa} {\rm Let $G$ be a finite subgroup of $SL(2,\k)$, and let $\k G$ be the group ring. $\k G$ acts on $\k[x,y]$ naturally. The skew group algebra $B=\k[x,y]\#\k G$ is a Koszul algebra. The projective dimension of the right graded $B$-module $B_0$ is two. Let $V=\k x\op\k y$ and $M=B_1=V\ot \k G$. Then $M\ot_{\k G}M\ot_{\k G}\cdots\ot_{\k G}M\cong V\ot\cdots\ot V \ot \k G$. With these notations, we have $T_{\k G}(M)\cong T(V)\#\k G$, and $B\cong T(V)\#\k G/(R)$ where the generating relation $R=\k (x\ot y-y\ot x)\ot \k G$. Since $G\subseteq SL(2,\k)$, the $\k G$-bimodule $R$ is isomorphic to the regular bimodule $\k G$. Hence any $\k G$-bimodule morphism $\theta:R\to \k G$ is determined by a central element $\lambda$ of $\k G$. By Theorem \ref{thm3}, we have $\mathbf{H}_\lambda =T(V)\#\k G/(x\ot y-y\ot x-\lambda)$ is a PBW deformation of $\k[x,y]\#\k G$. The filtered algebra $\mathbf{H}_\lambda$ is usually call a symplectic reflection algebra (see \cite{CBH,EG} when $\k$ is of characteristic zero, and \cite{N} when $\k$ is of positive characteristic).
}
\end{exa}

\begin{rem} {\rm (i) Note that in Theorem \ref{thm3}, the generating relations of the algebra $U$ do not involve the degree 1 part of the tensor algebra. As mentioned in the introduction, the setting of the theorem occurs quite common. Besides the example above, Theorem \ref{thm3} applies to many other interesting algebras. However, if the generating relations of $U$ involve degree 1 part as it did in \cite{BG,PP,SW,SW2,WW}, then the situation becomes more complicated. We can not prove the PBW theorem under the hypothesis on $B$ as in Theorem \ref{thm3}. We need more restrictions on $B$ as we will do in the next section.

(ii) If the projective dimension of $S_B$ is larger than two, then the method used in this section is not valid any more. To study PBW deformations of Koszul algebras with pdim$S_B\ge3$, or PBW deformations with generating relations involve degree 1 part, we will assume that $B$ is free as right $S$-module. Note that, in this case, $B_1\cong V\ot S$ for some vector space $V$. Thus, the left $S$-module action on $B_1$ yields a braiding $\Psi: S\ot V\to V\ot S$, which extends naturally to a braiding $\Psi_T:S\ot T(V)\to T(V)\ot S$ such that $(S,T(V),\Psi_T)$ forms an {\it entwining structure} (for the definition, see the next section). The graded algebra $B$ is a quotient algebra of a {\it generalized smash product} $T(V)\#S$. Therefore, we will consider the PBW deformations of a generalized smash product from a Koszul algebra together with a finite dimensional ring $S$. We leave details to the next section.}
\end{rem}

\section{PBW deformations of generalized smash products}\label{secc}

Let $B=\op_{i\ge 0}B_i$ be a graded algebra such that $\dim B_i<\infty$ for all $i\ge0$. Write $S$ for $B_0$, and $M$ for $B_1$. Assume that $B$ is generated in degree $1$ over $S$, that is, $B\cong T_S(M)/(\overline{R})$ for some finite dimensional graded sub-$S$-bimodule $\overline{R}\subseteq M\ot_SM\op M\ot_SM\ot_SM\op\cdots$. Assume further that $M$ is free as a right $S$-module. Then as a right $S$-module $M\cong V\ot S$ for some finite dimensional vector space $V$. Example \ref{ex1} below shows that the left $S$-module action on $M$ induces a {\it braiding} $\Psi:S\ot V\to V\ot S$. We may extend the braiding $\Psi$ naturally to  $\Psi_T:S\ot T(V)\to T(V)\ot S$ so that the triple $(S,T(V),\Psi_T)$ is an {\it entwining structure} (for the definition, see the next subsection). Then as a graded algebra $T_S(M)$ is isomorphic to $T(V)\#S$, the {\it smash product} of $T(V)$ with $S$ (cf. Lemma \ref{lem1}). Therefore $B$ is a quotient algebra of the smash product $T(V)\#S$, that is, there is a graded ideal $\overline{I}$ such that $B\cong (T(V)\#S)/\overline{I}$. Now assume that $B$ is a Koszul algebra. In many interesting situations (for instance, the graded algebras considered in \cite{EG,SW,SW2,N,NW}), the ideal $\overline{I}$ is generated by $R\ot S$ for some subspace $R\subset V\ot V$. Moreover, $A:=T(V)/(R)$ is a classical Koszul algebra and the braiding $\Psi_T$ induces a braiding $\overline{\Psi}:S\ot A\to A\ot S$ so that $(S,A,\overline{\Psi})$ is an entwining structure (see the next subsection). In this situation, the PBW deformations on $B$ are equivalent to the PBW deformations on $A\#S$. Hence, in this section, we only focus on the PBW deformations of the smash product algebras. We will first construct projective bimodule resolutions for $A\#S$ by using the similar techniques developed in \cite{SW}, and then generalize the PBW theorems in \cite{SW,SW2,WW} to our setting.

\subsection{Generalized smash products}\label{sec1}

Let $S$ be a finite dimensional algebra, $V$ a vector space. We call a linear map $$\Psi:S\ot V\longrightarrow V\ot S$$ a {\it braiding} if the following diagrams commute:
\begin{center}
$\xymatrix{
  S\ot S\ot V \ar[d]_{\mu_S\ot 1} \ar[r]^{1\ot\Psi} & S\ot V\ot S  \ar[r]^{\Psi\ot 1} & V\ot S\ot S\ar[d]^{1\ot \mu_S} \\
  S\ot V \ar[rr]^{\Psi} &  & V\ot S  }$ \hspace{10mm}
$\xymatrix{
  V \ar[d]_{\iota} \ar[dr]^{\gamma}        \\
  S\ot V\ar[r]^{\Psi}  & V\ot S ,           }$
\end{center}
where $\mu_S$ in the first diagram is the multiplication map of $S$, and the linear maps $\iota$ and $\gamma$ are defined by $\iota(v)=1\ot v$ and $\gamma(v)=v\ot 1$ for all $v\in V$.

\begin{exa}
{\rm (i) The twisting map $\tau: S\ot V\longrightarrow V\ot S$ is a braiding.

(ii) Let $H$ be a finite dimensional Hopf algebra, $M$ a left $H$-module. The linear map $\Psi: H\ot M\longrightarrow M\ot H$ by $\Psi(h\ot m)=h_{(1)}m\ot h_{(2)}$, for all $h\in H$ and $m\in M$, is a braiding, where we use Sweedler's notation \cite{S}.}
\end{exa}

\begin{exa}\label{ex1} {\rm Let $M$ be an $S$-bimodule which is free as a right $S$-module. As a right $S$-module, $M\cong V\ot S$ for some vector space $V$. We fix an isomorphism $\phi:M\to V\ot S$, through which the left $S$-module structure on $M$ induces a left $S$-module structure on $V\ot S$ so that $V\ot S$ is an $S$-bimodule. Define a linear map $\Psi: S\ot V\longrightarrow V\ot S$ by $\Psi(s\ot v)=s(v\ot 1)$ for all $s\in S$ and $v\in V$. Then $\Psi$ is a braiding. }
\end{exa}

Let $A$ be an algebra, $\Psi:S\ot A\longrightarrow A\ot S$ be a braiding. We call the triple $(S,A,\Psi)$ an {\it entwining structure} if the following diagram commute:
\begin{center}
$\xymatrix{
  S\ot A\ot A \ar[d]_{1\ot \mu_A} \ar[r]^{\Psi\ot1} & A\ot S\ot A  \ar[r]^{1\ot\Psi} & A\ot A\ot S\ar[d]^{\mu_A\ot 1} \\
  S\ot A \ar[rr]^{\Psi} &  & A\ot S  }$ \hspace{10mm}
$\xymatrix{
  S \ar[d]_{\nu} \ar[dr]^{\kappa}        \\
  S\ot A\ar[r]^{\Psi}  & A\ot S ,           }$
\end{center}
where $\mu_A$ is the multiplication map of $A$, and $\nu(s)=s\ot 1$ and $\kappa(s)=1\ot s$ for all $s\in S$.

\begin{note} {\rm Let $\Psi:S\ot V\longrightarrow V\ot S$ be a braiding. We use the following classical notation (cf. \cite[Sec. 2.3]{CMZ})  $$\Psi(v\ot s)=s^\Psi\ot v_\Psi$$ for $v\in V$ and $s\in S$, where the summation on the right hand side is understood to be omitted.}
\end{note}

Let $(S,A,\Psi)$ be an entwining structure. One may define a new algebra $A\#S$ \cite[Sec. 2.3, Theorem 7]{CMZ}. As a vector space, $A\#S=A\ot S$. The product is given by $$(a\ot s)(b\ot t)=a\ b^\Psi\ot s_\Psi\ t$$
for $a,b\in A$ and $s,t\in S$. The algebra $A\#S$ is called the {\it smash product} of $S$ and $A$ \cite[Sec. 2.3]{CMZ}.

Let $I$ be a two-sided idea of $A$ such that $\Psi(S\ot I)\subseteq I\ot S$. In this case,  $I\ot S$ is  a two-sided ideal of $A\#S$. The braiding $\Psi$ induces a new braiding $\overline{\Psi}:S\ot A/I\longrightarrow A/I\ot S$ so that $(S,A/I,\overline{\Psi})$ is an entwining structure. Hence we have the smash product $A/I\#S$. Moreover, we have the algebra isomorphism:
\begin{equation}\label{eq1}
  A/I\#S\cong A\#S/(I\ot S).
\end{equation}

Let $K$ be a right $A$-module. $K\#S$ is a right $A\#S$-module whose underlying vector space is $K\ot S$ and the right $A\#S$-action on $K\# S$ is defined by $$(k\ot s)(a\ot t)=k\ a^\Psi\# s_\Psi\ t$$ for $k\in K$, $a\in A$ and $s,t\in S$. Let Mod-$A$ be the category of right $A$-modules and Mod-$A\#S$ the category of right $A\#S$-modules. Then we have an exact functor
\begin{equation}\label{eq3}
  -\#S: \text{Mod-}A\longrightarrow \text{Mod-}A\#S,
\end{equation}
preserving projective modules.

Let $M$ be an $S$-bimodule such that it is free as a right $S$-module. Then there is a vector space $V$ and an isomorphism of right $S$-modules $\phi:M\longrightarrow V\ot S$. As we see in Example \ref{ex1}, the left $S$-action on $M$ induces a braiding $\Psi:S\ot V\longrightarrow V\ot S$. Let $T(V)$ be the tensor algebra of $V$. The braiding $\Psi$ induces a braiding $$\Psi_T:S\ot T(V)\longrightarrow T(V)\ot S$$ in the obvious way, that is, $\Psi_T(s\ot 1)=1\ot s$ for $s\in S$ and $\Psi_T$ acts on $S\ot V^{\ot n}$ ($n\ge1$) through the composition:
{\small \begin{equation*}
\xymatrix{
     S\ot V^{\ot n} \ar[rr]^{\Psi\ot 1^{\ot n-1}} && V\ot S\ot V^{\ot n-1} \ar[rr]^{\qquad 1\ot\Psi\ot 1^{\ot n-2}} &&\cdots\ar[r] &V^{\ot n-1}\ot S\ot V \ar[rr]^{1^{\ot n-1}\ot \Psi} &&V^{\ot n}\ot S. }
\end{equation*}}
Now it is straightforward to check that the triple $(S,T(V),\Psi_T)$ is an entwining structure. Therefore, we can form the smash product $T(V)\#S$. As before, we let $T_S(M)$ stand for the tensor algebra over $S$.

\begin{lem}\label{lem1} With the notations as above, there is an isomorphism of algebras $T_S(M)\cong T(V)\#S$.
\end{lem}
\proof Note that $V\ot S$ is  an $S$-bimodule, where the left $S$-module structure is induced by the left $S$-action on $M$ via the isomorphism $\phi$. Thus, $\phi$ is actually an isomorphism of $S$-bimodule. Now we may define a linear map $\Phi:T_S(M)\longrightarrow T(V)\#S$ as follows: $\Phi(s)=s$ for $s\in S$ and the restriction of $\Phi$ to $M\ot_S \cdots\ot_SM$ is the composition:
\begin{equation*}
\xymatrix{
     M\ot_S \cdots\ot_S M \ar[rr]^{\phi\ot_S\cdots\ot_S\phi\qquad} && (V\ot S)\ot_S\cdots\ot_S(V\ot S) \ar[r]^{\qquad\cong} &V\ot\cdots\ot V\ot S. }
\end{equation*}
Now it is a routine check that $\Phi$ is an isomorphism of algebras. \qed

\subsection{Generalized smash products and Koszul algebras}

Let $B$ be a Koszul algebra. Set $S=B_0$ and $M=B_1$. Then $B$ is a quadratic algebra in the sense that $B\cong T_S(M)/(\overline{R})$ where $\overline{R}$ is a sub-$S$-bimodule of $M\ot_SM$ and $(\overline{R})$ is the two-sided ideal of $T_S(M)$ generated by $\overline{R}$. If $M$ is free as a right $S$-module, then $T_S(M)\cong T(V)\#S$ for some vector $V$ and some braiding $\Psi:S\ot V\longrightarrow V\ot S$. Hence, in this case, $B$ is a quotient algebra of the smash product $T(V)\#S$.

Given a finite dimensional vector space and a bijective braiding $\Psi:S\ot V\longrightarrow V\ot S$, we get an entwining structure $(S,T(V),\Psi_T)$. Let $R\subseteq V\ot V$ be a subspace. Consider the right free $S$-module $R\ot S$. If $\Psi_T(S\ot R)\subseteq R\ot S$, then $R\ot S$ is a sub-$S$-bimodule of $T(V)\# S$. Let $I=(R)$ be the two-sided ideal of $T(V)$ generated by $R$. One can easily verify that $\Psi_T(S\ot I)\subseteq I\ot S$, so that $I\ot S$ is a two-sided idea of $T(V)\#S$. Let $A=T(V)/I$. We have an entwining structure $(S,A,\overline{\Psi})$ (cf. Subsection \ref{sec1}). Let $\overline{I}$ be the two-sided ideal of $T(V)\#S$ generated by $R\ot S$. One sees that $I\ot S=\overline{I}$. Therefore, we have the following isomorphisms of algebras:
\begin{equation}\label{eq2}
  T(V)\#S/\overline{I}\cong T(V)\#S/(I\ot S)\cong A\#S.
\end{equation}

\begin{prop}\label{prop1} With the notations as above. If $A$ is a classical Koszul algebra, then $A\#S$ is a Koszul algebra.
\end{prop}
\proof Since $\Psi$ is bijective, the braiding $\overline{\Psi}$ is also bijective. Hence $A\#S$ is free on both sides as an $S$-module. Since $A$ is classically Koszul, we have the Koszul resolution
\begin{equation}\label{eq4}
  \cdots\longrightarrow K_n\ot A\overset{\partial^{-n}}\longrightarrow\cdots\longrightarrow R\ot A\overset{\partial^{-2}}\longrightarrow V\ot A\overset{\partial^{-1}}\longrightarrow A\longrightarrow \k_A\longrightarrow0,
\end{equation}
where $K_n=\bigcap_{i=0}^{n-2}V^{\ot i}\ot R\ot V^{\ot n-i-2}\subseteq V^{\ot n}$ for $n\ge2$, and the differential $\partial^{-n}$ acts on pure tensors as follows: $\partial^{-n}(v_1\ot\cdots\ot v_n\ot a)=v_1\ot\cdots\ot v_{n-1}\ot v_na$ for $v_1,\dots,v_n\in V$ and $a\in A$. Applying the functor $-\#S$ (cf. Subsection \ref{sec1}) on the Koszul resolution above, we obtain an exact complex:
\begin{equation}\label{eq7}
  \cdots\longrightarrow K_n\ot A\#S\overset{\partial^{-n}\ot 1}\longrightarrow\cdots\longrightarrow R\ot A\#S\overset{\partial^{-2}\ot 1}\longrightarrow V\ot A\#S\overset{\partial^{-1}\ot 1}\longrightarrow A\#S\longrightarrow S\longrightarrow0.
\end{equation}
It is clear that $K_n\ot A\#S$, viewed as a graded module, is generated in degree $n$. Hence $A\#S$ is a Koszul algebra. \qed

There are three natural projective bimodule resolutions for the Koszul algebra $A\#S$. In the next three subsections, we will construct these resolutions and make a comparison of the resolutions.

\subsection{Projective bimodule resolution of $A\#S$ arising from the Koszul resolution of $A$}\label{sec3}

In this subsection, $S$ is an algebra, $V$ is a finite dimensional vector space, and $\Psi:S\ot V\longrightarrow V\ot S$ is a bijective braiding. We keep the same notations as in the previous subsection. Let $A=T(V)/(R)$ be a classical Koszul algebra. We assume that $\Psi_T(S\ot R)\subseteq R\ot S$. By Proposition \ref{prop1}, the smash product $A\#S$ is a Koszul algebra. We will construct a free $A\#S$-bimodule resolutions for  $A\#S$ from the Koszul resolution of $A$.

Since $A$ is classically Koszul, we have Koszul bimodule resolution of ${}_AA_A$:
\begin{equation}\label{eq8}
\cdots\longrightarrow A\ot K_n\ot A\overset{d^{-n}}\longrightarrow\cdots\longrightarrow A\ot R\ot A\overset{d^{-2} }\longrightarrow A\ot V\ot A\overset{d^{-1} }\longrightarrow A\ot A\overset{d^0}\longrightarrow A\longrightarrow0,
\end{equation}
where the maps $d^{-n}$ ($n\ge1$) are given  by
\begin{eqnarray*}
  &&d^{-n}(a\ot v_1\ot\cdots\ot v_n\ot b) \\
   &=&  av_1\ot v_2\ot\cdots\ot v_n\ot b+(-1)^na\ot v_1\ot\cdots\ot v_{n-1}\ot v_nb
\end{eqnarray*}
for $a,b\in A$, $v_1,\dots,v_n\in V$, and $d^0(a\ot b)=ab$.
Applying the functor $-\#S$ on the complex (\ref{eq8}), we obtain an exact sequence
{\small\begin{equation}\label{eq9}
\cdots\longrightarrow A\ot K_n\ot A\#S\overset{d^{-n}\ot 1}\longrightarrow\cdots\longrightarrow A\ot R\ot A\#S\overset{d^{-2}\ot 1}\longrightarrow A\ot V\ot A\#S\overset{d^{-1}\ot 1}\longrightarrow A\ot A\#S\overset{d^0\ot 1}\longrightarrow A\#S\longrightarrow0.
\end{equation}}
Each component in the exact sequence above is free as a right $A\#S$-module. The differentials are obviously compatible with the right $A\#S$-module structures. For $n\ge0$, we may define a left $A\#S$-action on $A\ot K_n\ot A\#S$ so that it is an $A\#S$-bimodule:
for $a,b\in A$, $x\in K_n$ and $z\in A\#S$, $$(a\ot s)(b\ot x\ot z)=ab^{\overline{\Psi}}\ot x^{\Psi_T}\ot s_{\overline{\Psi}\Psi_T}\ z.$$

\begin{lem}\label{lem3} The exact sequence (\ref{eq9}) is an exact complex of $A\#S$-bimodules.
\end{lem}
\proof We only need to verify that the differentials are left $S$-linear. Since any element $x\in K_n$ is linear combinations of pure tensors $v_1\ot\cdots\ot v_n$ with $v_1,\dots,v_n\in V$, it suffices to compute the following: for $s,s'\in S$, $a,b\in A$,
\begin{eqnarray*}
  &&(d^{-n}\ot1)[s(a\ot v_1\ot\cdots\ot v_n\ot b\ot s')]\\
  &=&(d^{-n}\ot1)[a^{\overline{\Psi}}\ot v_1^{\Psi}\ot x^{\Psi_T}\ot v_n^{\Psi}\ot b^{\overline{\Psi}}\ot s_{\overline{\Psi}\Psi\Psi_T\Psi}\ s'] \\
  &=& a^{\overline{\Psi}} v_1^{\Psi}\ot x^{\Psi_T}\ot v_n^{\Psi}\ot b^{\overline{\Psi}}\ot s_{\overline{\Psi}\Psi\Psi_T\Psi\overline{\Psi}}\ s'+(-1)^na^{\overline{\Psi}}\ot v_1^{\Psi}\ot x^{\Psi_T}\ot v_n^{\Psi} b^{\overline{\Psi}}\ot s_{\overline{\Psi}\Psi\Psi_T\Psi\overline{\Psi}}\ s',\\
  &&s[(d^{-n}\ot1)(a\ot v_1\ot\cdots\ot v_n\ot b\ot s')] \\
  &=& s(av_1\ot x\ot v_n\ot b\ot s'+(-1)^na\ot v_1\ot x\ot v_nb\ot s')\\
  &=&a^{\overline{\Psi}}v_1^{\Psi}\ot x^{\Psi_T}\ot v_n^\Psi b^{\overline{\Psi}}\ot s_{\overline{\Psi}\Psi\Psi_T\Psi\overline{\Psi}}\ s'+(-1)^na^{\overline{\Psi}}\ot v_1^{\Psi}\ot x^{\Psi_T}\ot v_n^{\Psi} b^{\overline{\Psi}}\ot s_{\overline{\Psi}\Psi\Psi_T\Psi\overline{\Psi}}\ s',
\end{eqnarray*}
where $x=v_2\ot\cdots\ot v_{n-1}\in T(V)$. Hence the left $S$-action is compatible with the morphism $d^{-n}\ot 1$. \qed

We next construct a free $A\#S$-bimodule resolution of $A\#S$ through the exact sequence (\ref{eq9}). Consider the bar resolution of the left $S$-module $K_n\ot S$:
\begin{equation}\label{eq10}
  \cdots\longrightarrow S\ot S^{\ot n}\ot K_n\ot S\longrightarrow\cdots\longrightarrow  S\ot K_n\ot S\longrightarrow K_n\ot S\longrightarrow0.
\end{equation}
Since $K_n\ot S$ is free as a right $S$-module, the above sequence is indeed a free $S$-bimodule resolution of $K_n\ot S$. Note that the sequence (\ref{eq10}) is split exact as a complex of right $S$-modules. Applying the functors $A\#S\ot_S-$ and $-\ot_S A\#S$ we obtain a free $A\#S$-bimodule resolution of $A\ot K_n\ot A\#S$:
\begin{equation}\label{eq11}
  \cdots\longrightarrow A\#S\ot S^{\ot m}\ot K_n\ot A\#S\overset{\partial^{-m}_n}\longrightarrow\cdots\overset{\partial^{-1}_n}\longrightarrow A\#S\ot K_n\ot A\#S\overset{\partial^0_n}\longrightarrow A\ot K_n\ot A\#S\longrightarrow0,
\end{equation}
where the differential $\partial^{-m}$ acts on pure tensors by
\begin{eqnarray*}
 \partial^{-m}_n(a\ot s_0\ot s_1\ot\cdots\ot s_m\ot x)&=&\sum_{i=0}^{m-1}(-1)^ia\ot s_0\ot\cdots\ot s_is_{i+1}\ot\cdots\ot s_m\ot x
\\ && +(-1)^ma\ot s_0\ot\cdots\ot s_{m-1}\ot s_m x,
\end{eqnarray*}
for $a\in A$, $s_0,\dots,s_m\in S$ and $x\in K_n\ot A\#S$.

For each $n\ge1$ and $m\ge0$, let $$\vartheta^{-n}_m:A\#S\ot S^{\ot m}\ot K_n\ot A\#S\longrightarrow A\#S\ot S^{\ot m}\ot K_{n-1}\ot A\#S$$ be the restriction to $A\#S\ot S^{\ot m}\ot K_n\ot A\#S$ of the $A\#S$-bimodule morphism:
 $$g^{-n}_m:A\#S\ot S^{\ot m}\ot V^{\ot n}\ot A\#S\longrightarrow A\#S\ot S^{\ot m}\ot V^{\ot n-1}\ot A\#S,$$ defined by
\begin{eqnarray*}
 &&g^{-n}_m(a\ot s_0\ot\cdots\ot s_m\ot v_1\ot\cdots\ot v_n\ot z)\\
 &=&av_1^{\Psi\cdots\Psi}\ot s_{0\Psi}\ot\cdots\ot s_{m\Psi}\ot v_2\ot\cdots\ot v_n\ot z\\
 &&+(-1)^na\ot s_0\ot\cdots\ot s_m\ot v_1\ot\cdots\ot v_{n-1}\ot v_nz,
\end{eqnarray*}
where $a\in A$, $s_0,\dots,s_m\in S$, $v_1,\dots,v_n\in V$ and $z\in A\#S$.
We have the following lemma.

\begin{lem}\label{lem2} {\rm(i)} For all $m,n\ge1$, $\vartheta^{-n+1}_m\partial^{-m}_n=\partial^{-m}_n\vartheta^{-n}_m$, and $(d^{-n}\ot1)\partial^0_n=\partial^0_n\vartheta^{-n}_0$;

{\rm(ii)} for all $n\ge1$ and $m\ge0$, $\vartheta^{-n+1}_m\vartheta^{-n}_m=0$.
\end{lem}

By Lemma \ref{lem2}, we obtain the following double complex where  the superscripts and subscripts of the differentials are omitted:
{\tiny
$$\xymatrix{
 &0&0&0\\
 \cdots\ar[r]& A\#S\ot R\ot A\#S \ar[u] \ar[r]^{\vartheta} & A\#S\ot V\ot A\#S \ar[u] \ar[r]^{\vartheta} & A\#S\ot A\#S\ar[u]\ar[r]&0 \\
 \cdots\ar[r]& A\#S\ot  S\ot R\ot A\#S \ar[u]_{\partial} \ar[r]^{\vartheta} & A\#S\ot S\ot V\ot A\#S \ar[u]^{-\partial} \ar[r]^{\vartheta} & A\#S\ot S\ot A\#S \ar[u]^{\partial}\ar[r]&0 \\
   \cdots\ar[r]& A\#S\ot S\ot S\ot R\ot A\#S \ar[u]_{\partial} \ar[r]^{\vartheta} & A\#S\ot S\ot S\ot V\ot A\#S \ar[u]^{-\partial} \ar[r]^{\vartheta} & A\# S\ot S\ot S\ot A\#S \ar[u]^{\partial}\ar[r]&0\\
    &\vdots\ar[u]&\vdots\ar[u]&\vdots\ar[u]&  }
$$}

Let $P^\bullet$ be the total complex of the double complex above. Taking the $0^{th}$ cohomology of the vertical complexes of the double complex yields the following complex:
$$\cdots\longrightarrow A\ot K_n\ot A\#S\overset{d^{-n}\ot 1}\longrightarrow\cdots\longrightarrow A\ot R\ot A\#S\overset{d^{-2}\ot 1}\longrightarrow A\ot V\ot A\#S\overset{d^{-1}\ot 1}\longrightarrow A\ot A\#S\longrightarrow0.$$
By Lemma \ref{lem3}, this complex is exact except at the final position, and the $0^{th}$ cohomology of this complex is $A\#S$. Hence, by taking the classical spectral sequence of $P^\bullet$ we see that $P^\bullet$ is indeed a free $A\#S$-bimodule resolution of $A\#S$.

In summary, we obtain the following result, which may be viewed as a generalization of \cite[Theorem 4.3]{SW} and of \cite[Theorem 2.10]{WW}. Note that in \cite[Theorem 4.3]{SW} (resp. \cite[Theorem 2.10]{WW}), the Koszul resolution is replaced by a free resolution of $A$ which also has a module structure over the group ring (resp. the Hopf algebra) under consideration. Projective resolutions for smash product algebras over a Hopf algebra were also obtained in \cite{Ne}.

\begin{thm}\label{thm} $(P^\bullet,\delta)$ is a free $A\#S$-bimodule resolution of $A\#S$. Explicitly, for $n\ge0$, $$P^{-n}=\bigoplus_{i=0}^n A\#S\ot S^{\ot i}\ot K_{n-i}\ot A\#S$$ and the differential $\delta^{-n}$ $(n\ge1)$ acts on $A\#S\ot S^{\ot i}\ot K_{n-i}\ot A\#S$ by $$\delta^{-n}(A\#S\ot S^{\ot i}\ot K_{n-i}\ot A\#S)=\left((-1)^i\partial^{-i}_{n-i}+\vartheta^{n-1}_i\right)(A\#S\ot S^{\ot i}\ot K_{n-i}\ot A\#S).$$
\end{thm}

\subsection{Projective resolution of $A\#S$ arising from the bar resolution of $A$}\label{subsecc1}
Keep the notions as in the previous subsection.
We construct another free resolution of $A\#S$. In this subsection, set $J=\op_{i\ge1}A_i$. We begin with the bar resolution of $A$:
$$\cdots\longrightarrow A\ot J^{\ot n}\ot A\overset{\overline{d}^{-n}}\longrightarrow\cdots\overset{\overline{d}^{-2}}\longrightarrow A\ot J\ot A\overset{\overline{d}^{-1}}\longrightarrow A\ot A\overset{\overline{d}^0}\longrightarrow A\longrightarrow0.$$ Applying the functor $-\#S$ to the sequence above, we obtain an exact sequence:
\begin{equation}\label{eq12}
\cdots\longrightarrow A\ot J^{\ot n}\ot A\#S\overset{\overline{d}^{-n}\ot1}\longrightarrow\cdots\overset{\overline{d}^{-2}\ot1}\longrightarrow A\ot J\ot A\#S\overset{\overline{d}^{-1}\ot1}\longrightarrow A\ot A\#S\overset{\overline{d}^0\ot1}\longrightarrow A\#S\longrightarrow0,
\end{equation}
where the maps $\overline{d}^{-n}\ot1$ ($n\ge0$) acts on pure tensors by
\begin{eqnarray*}
  &&d^{-n}(a_0\ot a_1\ot\cdots\ot a_n\ot a_{n+1}\ot s) \\
   &=&  \sum_{i=0}^{n}(-1)^ia_0\ot 1_1\ot\cdots\ot a_ia_{i+1}\cdots\ot a_n\ot a_{n+1}\ot s
\end{eqnarray*}
for $a_0,\dots,a_{n+1}\in A$.

In the above sequence, each component is obviously a right $A\#S$-module, and each morphism is a right $A\#S$-module morphism. Thus, we may define a left $A\#S$-module structure on $A\ot J^{\ot n}\ot A\#S$ ($n\ge0$) using the entwining structure of $(A,S)$, that is, for $a,a_{n+1}\in A$, $a_1,\dots,a_n\in J$ and $s,s'\in S$, $$(a\ot s)(a_0\ot\cdots \ot a_{n+1}\ot s')=aa_0^{\overline{\Psi}\cdots\overline{\Psi}}\ot \cdots\ot a_{n+1}^{\overline{\Psi}}\ot s_{\overline{\Psi}\cdots\overline{\Psi}_T}\ s'.$$ Now it is straightforward to check that $A\ot J^{\ot n}\ot A\#S$ is an $A\#S$-bimodule and the morphisms in the seqence (\ref{eq12}) are $A\#S$-bimodule morphisms.

Similar to the previous case, we can construct a free resolution for $A\#S$ from the sequence (\ref{eq12}). For $n\ge1$, we consider the bar resolution of the left $S$-module $J^{\ot n}\ot S$:
\begin{equation}\label{eq13}
  \cdots\longrightarrow S\ot S^{\ot m}\ot J^{\ot n}\ot S\longrightarrow\cdots\longrightarrow  S\ot J^{\ot n}\ot S\longrightarrow J^{\ot n}\ot S\longrightarrow0,
\end{equation}
which is in fact a free resolution of the $S$-bimodule $J^{\ot n}\ot S$. Applying the functors $A\#S\ot _S-$ and $-\ot_S A\#S$ on the sequence above simultaneously, we obtain a free resolution of the $A\#S$-bimodule $A\ot J^{\ot n}\ot A\#S$:
\begin{equation}\label{eq14}
  \cdots\longrightarrow A\#S\ot S^{\ot m}\ot J^{\ot n}\ot A\#S\overset{\overline{\partial}^{-m}_n}\longrightarrow\cdots\overset{\overline{\partial}^{-1}_n}\longrightarrow A\#S\ot J^{\ot n}\ot A\#S\overset{\overline{\partial}^0_n}\longrightarrow A\ot J^{\ot n}\ot A\#S\longrightarrow0,
\end{equation}
where the differential $\overline{\partial}^{-m}$ acts on pure tensors by
\begin{eqnarray*}
 \overline{\partial}^{-m}_n(a\ot s_0\ot s_1\ot\cdots\ot s_m\ot x)&=&\sum_{i=0}^{m-1}(-1)^ia\ot s_0\ot\cdots\ot s_is_{i+1}\ot\cdots\ot s_m\ot x
\\ && +(-1)^ma\ot s_0\ot\cdots\ot s_{m-1}\ot s_m x,
\end{eqnarray*}
for $a\in A$, $s_0,\dots,s_m\in S$ and $x\in J^{\ot n}\ot A\#S$.

For each $n\ge1$ and $m\ge0$, we define an $A\#S$-bimodule morphism $$\overline{\vartheta}^{-n}_m:A\#S\ot S^{\ot m}\ot J^{\ot n}\ot A\#S\longrightarrow A\#S\ot S^{\ot m}\ot J^{\ot n-1}\ot A\#S$$ by
\begin{eqnarray*}
 &&\overline{\vartheta}^{-n}_m(a\ot s_0\ot\cdots\ot s_m\ot a_1\ot\cdots\ot a_n\ot z)\\
 &=&aa_0^{\overline{\Psi}\cdots\overline{\Psi}}\ot s_{0\overline{\Psi}}\ot\cdots\ot s_{m\overline{\Psi}}\ot a_2\ot\cdots\ot a_n\ot z\\
 &&+\sum_{i=1}^{n-1}(-1)^ia\ot s_0\ot\cdots\ot s_m\ot a_1\ot\cdots\ot a_ia_{i+1}\ot\cdots\ot a_{n}\ot z\\
 &&+(-1)^na\ot s_0\ot\cdots\ot s_m\ot a_1\ot\cdots\ot a_{n-1}\ot a_nz,
\end{eqnarray*}
where $a\in A$, $s_0,\dots,s_m\in S$, $a_1,\dots,a_n\in J$ and $z\in A\#S$.

Similar to Lemma \ref{lem2} and Theorem \ref{thm}, we obtain the following results,  generalizing partially \cite[Theorem 4.3]{SW} and \cite[Theorem 2.10]{WW}.
\begin{lem}\label{lem4} {\rm(i)} For all $m,n\ge1$, $\overline{\vartheta}^{-n+1}_m\overline{\partial}^{-m}_n=\overline{\partial}^{-m}_n\overline{\vartheta}^{-n}_m$, and $(d^{-n}\ot1)\overline{\partial}^0_n=\overline{\partial}^0_n\overline{\vartheta}^{-n}_0$;

{\rm(ii)} for all $n\ge1$ and $m\ge0$, $\overline{\vartheta}^{-n+1}_m\overline{\vartheta}^{-n}_m=0$.
\end{lem}

Lemma \ref{lem4} leads to the following double complex where the superscripts and subscripts of the differentials are omitted:
{\tiny
$$\xymatrix{
 &0&0&0\\
 \cdots\ar[r]& A\#S\ot J\ot J\ot A\#S \ar[u] \ar[r]^{\overline{\vartheta}} & A\#S\ot J\ot A\#S \ar[u] \ar[r]^{\overline{\vartheta}} & A\#S\ot A\#S\ar[u]\ar[r]&0 \\
 \cdots\ar[r]& A\#S\ot  S\ot J\ot J\ot A\#S \ar[u]_{\overline{\partial}} \ar[r]^{\overline{\vartheta}} & A\#S\ot S\ot J\ot A\#S \ar[u]^{-\overline{\partial}} \ar[r]^{\overline{\vartheta}} & A\#S\ot S\ot A\#S \ar[u]^{\overline{\partial}}\ar[r]&0 \\
   \cdots\ar[r]& A\#S\ot S\ot S\ot J\ot J\ot A\#S \ar[u]_{\overline{\partial}} \ar[r]^{\overline{\vartheta}} & A\#S\ot S\ot S\ot J\ot A\#S \ar[u]^{-\overline{\partial}} \ar[r]^{\overline{\vartheta}} & A\# S\ot S\ot S\ot A\#S \ar[u]^{\overline{\partial}}\ar[r]&0\\
    &\vdots\ar[u]&\vdots\ar[u]&\vdots\ar[u]&  }
$$}
\begin{prop}\label{prop3} Let $(\overline{P}^\bullet,\overline{\delta})$ be the total complex of the double complex above. Then it is a free $A\#S$-bimodule resolution of $A\#S$. Explicitly, for $n\ge0$, $$\overline{P}^{-n}=\bigoplus_{i=0}^n A\#S\ot S^{\ot i}\ot J^{\ot n-i}\ot A\#S$$ and the differential $\overline{\delta}^{-n}$ $(n\ge1)$ acts on $A\#S\ot S^{\ot i}\ot J^{\ot n-i}\ot A\#S$ by $$\delta^{-n}(A\#S\ot S^{\ot i}\ot J^{\ot n-i}\ot A\#S)=\left((-1)^i\overline{\partial}^{-i}_{n-i}+\overline{\vartheta}^{n-1}_i\right)(A\#S\ot S^{\ot i}\ot J^{\ot n-i}\ot A\#S).$$
\end{prop}

\subsection{Comparing projective resolutions of $A\#S$}\label{subsecc2}

As a matter of fact, the free resolution of $A\#S$ constructed in Theorem \ref{thm} is a subcomplex of the free resolution constructed in Proposition \ref{prop3}. We can construct an inclusion map from first resolution to the second one. Recall that $A=T(V)/(R)$. We identify $A_1$ with $V$. Then $K_n\subseteq A_1^{\ot n}$. Hence, for $m,n\ge1$, the free module $A\#S\ot S^{\ot m}\ot K_n\ot A\#S$ is a submodule of $A\#S \ot S^{\ot m}\ot J^{\ot n}\ot A\#S$. Let $\zeta_{m,n}$ be the inclusion map and $\zeta_{0,0}$ be the identity map. Since the entwining structure $(A,S,\overline{\Psi})$ is induced from the braiding $\Psi:S\ot V\to V\ot S$,
we see that $\Psi(s\ot v)=\overline{\Psi}(s\ot v)$, for $v\in V=A_1$ and $s\in S$. Now one sees that the inclusion maps $\zeta_{m,n}$ is compatible with the differentials in the double complexes above. Therefore, we obtain an inclusion map: $\zeta:P^\bullet\to\overline{P}^\bullet$.

Summarizing the arguments above, we obtain the following:
\begin{prop} \label{prop4} We have an inclusion map of $A\#S$-bimodule complex  $$\zeta:(P^\bullet,\delta)\longrightarrow(\overline{P}^\bullet,\overline{\delta}).$$
\end{prop}

The third natural free bimodule resolution of $A\#S$ is the bar resolution. Write $B$ for the graded algebra $A\#S$. Recall from Subsection \ref{sec0} that the bar resolution of $B$ is denoted by $Q^\bullet$. Similar to \cite[Lemma 4.7]{SW}, we have the following result.

\begin{prop}\label{prop5} There are morphisms of cochain complexes $\psi:P^\bullet\to Q^\bullet$ and $\varphi:Q^\bullet\to P^\bullet$ such that $\varphi\psi$ is homotopic to $id_{P^\bullet}$ and $\psi\varphi$ is homotopic to $id_{Q^\bullet}$. Moreover, $\psi\varphi$ is the identity map when restricted to $B\ot K_n\ot B$ for all $n\ge0$.
\end{prop}
\proof The same proof of \cite[Lemma 4.7]{SW} applies in our situation. \qed

\subsection{PBW Deformations} \label{sec4}
Let $A$ and $S$ be the same as in Subsections \ref{sec3}--\ref{subsecc2}, and let $B=A\#S$. We view $S$ and $A$ as subalgebras of $B$ in the obvious way, and treat elements of $S$ and $A$ as  elements of $B$ without mentioning the inclusion maps.

Let $M=V\ot S$. By Lemma \ref{lem1}, there is an isomorphism of graded algebras $\Phi:T_S(M)\to T(V)\#S$. Let $\overline{R}=R\ot S$. We have already seen in Subsection \ref{sec1}, that $T_S(M)/(\overline{R})\cong A\#S$. Let $\phi:R\to V\ot S$ and $\theta:R\to S$ be two linear maps such that
\begin{equation}\label{eq20}
  \phi(r^{\Psi_T})s_{\Psi_T}=s\phi(r)\text{ and }\theta(r^{\Psi_T})s_{\Psi_T}=s\theta(r)
\end{equation}
for all $r\in R$ and $s\in S$. These two maps induce two $S$-bimodule morphism respectively: $\phi_B:R\ot S\to V\ot S$ and $\theta_B:R\ot S\to S$ given by  $\phi_B(r\ot s)=\phi(r)s$ and $\theta_B(r\ot s)=\theta(r)s$ for all $r\in R$ and $s\in S$. Denote by $U$ the algebra $T_S(M)/(r-\phi_B(r)-\theta_B(r):r\in R)\cong T(V)\#S/(r-\phi_B(r)-\theta_B(r):r\in R)$.

We now have the following result, which is a generaliztion of \cite[theorem 5.4]{SW} from the case that $S$ is a group ring, and of \cite[Theorem 0.4]{WW} from the case that $S$ is a Hopf algebra.

\begin{thm}\label{thm2} $U$ is a PBW deformation of $B=A\#S$ if and only if the following conditions are satisfied:

{\rm(i)} $(\phi_B\ot_S \id-\id\ot_S\phi_B)(\overline{R}\ot_S M\cap M\ot_S\overline{R})\subseteq \overline{R}$;

{\rm(ii)} $\phi_B(\phi_B\ot_S \id-\id\ot_S\phi_B)=-(\theta_B\ot_S \id-\id\ot_S \theta_B)$;

{\rm(iii)} $\theta_B(\id\ot_S\phi_B-\phi_B\ot_S \id)=0$,

where the maps in {\rm(ii)} and {\rm (iii)} are defined over $\overline{R}\ot_SM\cap M\ot_S\overline{R}$.
\end{thm}
\proof The necessary condition has been proved in \cite[Lemma 5.2]{SW}. For the sufficiency, we note that a PBW deformation of $B$ maybe obtained from a graded deformation $B_t$ by specializing at $t=1$ \cite{BG}. Theorem \ref{thm} says that the free (graded) $A\# S$-module $P^{-n}$ is generated by elements in degrees $0,1,\dots,n$, implying that the $n$th Hochschild cohomology $HH^n(B,B)$ of the graded algebra $B$ is concentrated in degrees not less than $-n$. In particular, the third Hochschild cohomology $HH^3(B,B)$ is concentrated in degrees not less than $-3$. Hence by Proposition \ref{prop6}, if we could construct a third level graded deformation of $B$, then it may lift to a graded deformation automatically. Now the rest of the proof is exactly the same as that of \cite[Theorem 5.4]{SW}. We only remark that the requirements (\ref{eq20}) insure that $\phi$ and $\theta$ are 2-cocycles in the cochain complex $\underline{\Hom}_{B^e}(P^\bullet,B)$, which is homotopically equivalent to the Hochschild cochain complex $CH^\bullet(B,B)$ (cf. Subsection \ref{subsec0}) via the morphisms constructed in Proposition \ref{prop4}. \qed

\section{Deformations of Koszul Artin-Schelter Gorenstein algebras of injective dimension 2}\label{sec5}

\subsection{Standard Koszul Artin-Schelter Gorenstein algebras}

Let $S$ be a finite dimensional algebra, $B$ a locally finite positively graded algebra such that $B_0=S$. We call $B$ an {\it Artin-Schelter Gorenstein} algebra \cite{MM} if
\begin{itemize}
\item[(i)] the graded injective dimension injdim$B_B$=injdim${}_BB=d<\infty$,
\item[(ii)] as a graded $B$-module$$\underline{\Ext}_B^i(S_B,B)\cong\left\{
\begin{array}{ll}
0, & \hbox{if $i\neq d$;} \\
D(S)(l), & \hbox{if $i=d$,}
\end{array}
\right.
$$ where $D(S)=\Hom(S,\k)$ is viewed as a left graded $B$-bimodule concentrated in degree zero,
\item[(iii)] the left version of (ii) holds.
\end{itemize}
The integer $l$ is usually called the {\it Gorenstein parameter} of $B$.
If, in addition, $B$ is a Koszul algebra and the Gorenstein parameter of $B$ equals to the graded injective dimension of $B$, then $B$ is called a {\it standard} Koszul Artin-Schelter Gorenstein algebra.

\begin{lem} \label{lem7} Let $B$ be a standard Koszul Artin-Schelter Gorenstein algebra of graded injective dimension $d$. Then
\begin{itemize}
\item[\rm(i)] $S$ is selfinjective;
\item[\rm(ii)] ${}_BD(S)$ (resp. $D(S)_B$) is a Koszul $B$-module with projective dimension pdim${}_BD(S)=d$ (resp. pdim$D(S)_B=d$).
\end{itemize}
\end{lem}
\proof (i) Applying the functor $\underline{\Hom}_{B}(-,B)$ to the Koszul resolution of $S_B$ (cf. Proposition \ref{prop7}), we obtain a sequence of left graded $B$-module
\begin{equation}\label{eq31}
  0\longrightarrow B\longrightarrow\cdots\overset{{\partial^{-d+1}}^*}\longrightarrow B\ot_S K_{d-1}^\vee\overset{{\partial^{-d}}^*}\longrightarrow B\ot_S K_d^\vee\overset{{\partial^{-d-1}}^*}\longrightarrow B\ot_S K_{d+1}^\vee\longrightarrow\cdots,
\end{equation}
where $K_i^\vee=\Hom_S(K_i,S)$ for all $i\ge1$. Since $B\ot_S K_{d-1}^\vee$ is concentrated in degrees not less than $-d+1$, we see the $-d$th component of the left graded $S$-module $\underline{\Ext}_B^d(S_B,B)$ is a submodule of $K_{d}^\vee$. By the Gorenstein condition of $B$, we obtain that the left $S$-mdoule $D(S)$ is isomorphic to the kernel of the restriction of ${\partial^{-d}}^*$ to $K^\vee_d$ (by ignoring the grading). Since $D(S)$ is an injective $S$-module, it is isomorphic to a direct summand of $K_d^\vee$ (by ignoring the grading). Hence the left $S$-module $D(S)$ is also projective since $K^\vee_d$ is a projective $S$-module by the Koszulity of $B$. Therefore $S$ is selfinjective.

(ii) Let $L_0=\ker({\partial^{-d}}^*|_{K_{d}^\vee})$. Since ${}_BD(S)(d)$ is isomorphic to $L_0$, we may choose a submodule $N$ of $K^\vee_d$ so that $K^\vee_d=L_0\op N_0$. Hence $B\ot_S K_d^\vee=B\ot_S L_0\bigoplus B\ot_SN_0$. Since $B\ot_S L_0$ is contained in the kernel of ${\partial^{-d}}^*$, there is an $S$-submodule $U$ of $B_1\ot_SN_0$ such that $\ker({\partial^{-d}}^*)_{-d+1}=B_1\ot_S L_0\op U$. Since $(-d+1)$th cohomology of the complex (\ref{eq31}) is zero, we see that there is a decomposition of the $S$-module $K_{d-1}^\vee=L_1\op N_1$ such that ${\partial^{-d+1}}^*(L_1)=B_1\ot_SL_0$ and ${\partial^{-d+1}}^*(N_1)=U$. Since $\underline{\Ext}_{B}^d(S_B,B)\cong DS(d)$, the sequence $B\ot_SL_1\overset{{\partial^{-d+1}}^*}\longrightarrow B\ot L_0\longrightarrow DS\longrightarrow0$ is exact. Inductively, we may find $S$-modules $L_i$ and $N_i$ ($2\leq i\leq d$) such that $K^\vee_{d-i}=L_i\op N_i$, ${\partial^{-d+i-1}}^*(L_i)\subseteq B_1\ot_S L_{i-1}$ and ${\partial^{-d+i-1}}^*(L_i)\subseteq B_1\ot_S N_i$, and the sequence $0\longrightarrow B\ot_S L_d\overset{{\partial^{-1}}^*}\longrightarrow B\ot_SL_{d-1}\longrightarrow\cdots\overset{{\partial^{-d+1}}^*}\longrightarrow B\ot_S L_0\longrightarrow D(S)(d)\longrightarrow0$ is exact. We claim that $L_i\neq0$ for all $0\leq i\leq d$. If $L_k=0$ for some $0\leq k\leq d$, than $L_i=0$ for all $i\ge k$ by the construction above. Hence the projective dimension of ${}_BD(S)$ is less than $d$. Since $S$ is selfinjective, any indecomposable projective left $S$-module is isomorphic to a direct summand of $D(S)$. If we write the graded $B$-module ${}_BS$ in a direct sum of indecomposable modules, than each component is isomorphic to a direct summand of ${}_BD(S)$. Hence the graded projective dimension of ${}_BS$ must be less than $d$, which contradicts to the Gorenstein conditions. Hence the claim follows. Therefore ${}_BD(S)$ is a Koszul $B$-module of projective dimension $d$. Symmetrically, we see that $D(S)_B$ is also a Koszul module of projective dimension $d$. \qed

Lemma \ref{lem7} implies the following result which generalizes \cite[Theorem 4.11]{DW} (also see \cite[Proposition 3.3]{MW}).

\begin{prop} \label{prop10} Let $B$ be a standard Koszul Artin-Schelter Gorenstein algebra of injective dimension $d$. Then $S_B$ and ${}_BS$ are of projective dimension $d$.
\end{prop}
\proof Let $P_1,\dots,P_k$ be the set of nonisomorphic indecomposable right projective $S$-modules. Since $S$ is selfinjective, $D(S)$ viewed as a right $S$-module is a direct sum of copies of $P_1,\dots,P_k$. Now Lemma \ref{lem7} implies that $S_B$ is of projective dimension $d$. Similarly, the projective dimension of ${}_BS$ is $d$. \qed

\subsection{Deformations of standard Koszul Artin-Schelter Gorenstein algebras of graded injective dimension 2}
Let $M$ be a finite dimensional $S$-bimodule, $R\subseteq M\ot_S M$ be a sub-$S$-bimodule, and let $B=T_S(M)/(R)$. Assume that $B$ is a standard Koszul Artin-Schelter Gorenstein algebra of graded injective dimension 2.

\begin{lem}\label{lem8} If $S$ is a basic algebra or $S$ is Frobenius, then there is an automorphism $\sigma$ of the algebra $S$ such that $R\cong {}_\sigma S_1$ as $S$-bimodules, where ${}_\sigma S_1$ is the $S$-bimodule with the right regular $S$-action and the left $S$-action twisted by $\sigma$.
\end{lem}
\proof By Lemma \ref{lem7}, $S$ is selfinjective. If $S$ is a basic algebra, then $S$ is Frobenius. By Proposition \ref{prop10}, $S_B$ is of projective dimension 2. Hence the Koszul resolution of $S_B$ reads as $0\longrightarrow R\ot_SB\longrightarrow M\ot_SB\longrightarrow B\longrightarrow S_B\longrightarrow0$. Since $B$ is a standard Koszul Artin-Schelter Gorenstein algebra, $\Hom_S(R_S,S)\cong D(S)$ as left $S$-modules. Let $\tau$ be the Nakayama automorphism of $S$ as $S$ is Frobenius. Then ${}_\tau S_1\cong D(S)$ as $S$-bimodules. Hence we have left $S$-module isomorphisms $\Hom_S(R_S,S)\cong \Hom_S(R_S,D({}_\tau S_1))\cong D(R_{\tau^{-1}})$. Therefore as a right $S$-module $R_{\tau^{-1}}\cong S$. Similarly, we have a left $S$-module isomorphism ${}_\tau R\cong S$. Hence the $S$-bimodule ${}_\tau R_{\tau^{-1}}$ is isomorphic to $S$ both as a left $S$-module and as a right $S$-module. Thus, there is an algebra automorphism $\phi$ of $S$ such that ${}_\tau R_{\tau^{-1}}\cong {}_\phi S_1$ as $S$-bimodules (cf. \cite[Lemma 2.9]{MM}). Now let $\sigma=\tau^{-1}\phi\tau^{-1}$. Then $R\cong {}_\sigma S_1$. \qed

Note that the automorphism $\sigma$ in Lemma \ref{lem8} is unique up to inner automorphism of $S$. Let $f:{}_\sigma S_1\overset{\cong}\longrightarrow R$ be the isomorphism in the lemma above. Set $r_0=f(1)$. By Theorem \ref{thm3}, we have the following result.

\begin{thm}\label{thm4} Let $S$ and $\sigma$ be the same as in Lemma \ref{lem8}. For any element $e\in S$ such that $se=e\sigma(s)$ for all $s\in S$, $U_e=T_S(M)/(r_0-e)$ is a PBW deformation of $B$.
\end{thm}

\vspace{5mm}

\subsection*{Acknowledgement} The authors are very grateful to the referee for his/her careful reading of the manuscript, and valuable comments and suggestions, which greatly enhance the readability and the interest of the article. The first author is supported by a grant from NSFC (No. 11171067), the second and the third authors are supported by FWO.

\vspace{5mm}

\bibliography{}

\end{document}